\documentclass[11pt,reqno]{amsart}
\usepackage[skip=5pt plus1pt, indent=20pt]{parskip}

\usepackage{palatino}
\usepackage{mathpazo}
\usepackage{amsmath,amssymb,amsxtra,color,calligra,mathrsfs,tcolorbox,stmaryrd}

\usepackage{xcolor} % Required for specifying custom colours
\colorlet{mdtRed}{red!50!black}
\colorlet{dblue}{blue!50!black}
\usepackage[colorlinks,pagebackref=true]{hyperref}
\hypersetup{linkcolor=dblue,citecolor=dblue,filecolor=dullmagenta,urlcolor=mdtRed}
\renewcommand*{\backref}[1]{}
\renewcommand*{\backrefalt}[4]{[{%
		\ifcase #1 Not cited.%
		\or $\uparrow$~#2.%
		\else $\uparrow$~#2.%
		\fi%
	}]}
\usepackage[all]{xy}
\usepackage{tikz,tikz-cd,tkz-graph,enumerate}
\usetikzlibrary{matrix,arrows,decorations.pathmorphing}

\usepackage[us,12hr]{datetime}
\usepackage{comment}

\DeclareMathOperator{\Spec}{{\rm Spec}}
\DeclareMathOperator{\Id}{{\rm Id}}
\DeclareMathOperator{\Hom}{{\rm Hom}}

\DeclareMathOperator{\Sch}{{\rm Sch}}

\DeclareMathOperator{\Vect}{{\bf Vect}}
\DeclareMathOperator{\PVect}{{\bf ParVect}}

\DeclareMathOperator{\Con}{{\rm Con}}

\DeclareMathOperator{\Mod}{{\mathfrak{Mod}}}
\DeclareMathOperator{\Coh}{{\mathfrak{Coh}}}

\DeclareMathOperator{\Filt}{{\bf Filt}}

\newcommand{\mf}[1]{\mathfrak{#1}}
\newcommand{\mc}[1]{\mathcal{#1}}
\newcommand{\ms}[1]{\mathscr{#1}}
\newcommand{\bb}[1]{\mathbb{#1}}

\newcommand{\doi}[1]{\href{https://doi.org/#1}{doi:#1}}
\numberwithin{equation}{subsection}

\newtheorem{theorem}[equation]{Theorem}
\newtheorem{corollary}[equation]{Corollary}
\newtheorem{lemma}[equation]{Lemma}
\newtheorem{proposition}[equation]{Proposition}

\theoremstyle{definition}
\newtheorem{definition}[equation]{Definition}
\newtheorem{remark}[equation]{Remark}

\theoremstyle{theorem}

\makeatletter
\newcommand\fnsymb[1]{\textsuperscript{\@fnsymbol{#1}}}
\newcommand\fnletter[1]{\lowercase{\textsuperscript{\@alph{#1}}}}

\makeatother

\usepackage{marvosym} % For Email/Letter icons. 
\renewcommand{\email}[2][1]{\thanks{\textit{Email address}#1: \href{mailto:#2}{#2}}}

\renewcommand{\address}[2][1]{\thanks{\textit{Address}#1: #2}} 

\begin{document}

\baselineskip=15.5pt 

\title[Equivariant Parabolic connections and stack of roots]{Equivariant Parabolic connections and stack of roots}

\author[S. Chakraborty]{Sujoy Chakraborty\fnsymb{1}}

\address[\fnsymb{1}]{Department of Mathematics, Indian Institute of Science Education and Research Tirupati, Andhra Pradesh - 517619, India.} 

\email[\fnsymb{1}]{sujoy.cmi@gmail.com}

\author[A. Paul]{Arjun Paul\fnsymb{2}}

\address[\fnsymb{2}]{Department of Mathematics and Statistics, 
	Indian Institute of Science Education and Research Kolkata, 
	Mohanpur - 741 246, Nadia, West Bengal, India.} 

\email[\fnsymb{2}]{arjun.paul@iiserkol.ac.in}

\thanks{Corresponding author: Arjun Paul}



\begin{abstract}
	Let $X$ be a smooth complex projective variety equipped with an action of a linear algebraic group $G$ over $\mathbb{C}$. Let $D$ be a reduced effective divisor on $X$ that is invariant under the $G$--action on $X$. Let $s_D$ be the canonical section of $\mathcal{O}_X(D)$ vanishing along $D$. Given a positive integer $r$, consider the stack $\mathfrak{X} := \mathfrak{X}_{(\mathcal{O}_X(D),\, s_D,\, r)}$ of $r$-th roots of $(\mathcal{O}_X, s_D)$ together with the natural morphism $\pi : \mathfrak{X} \to X$. Under the assumption that $G$ has no non-trivial characters, we show that the $G$--action on $X$ naturally lifts to a $G$--action on $\mathfrak{X}$ such that $\pi$ become $G$--equivariant, and the tautological invertible sheaf $\mathscr{M}$ on $\mathfrak{X}$ admits a linearization of this $G$--action. Finally, we define the notions of $G$--equivariant logarithmic connections on $\mathfrak{X}$ and $G$--equivariant parabolic connections on $X$ with rational parabolic weights along $D$, and establish an equivalence between the category of $G$--equivariant logarithmic connections on $\mathfrak{X}$ and the category of $G$--equivariant parabolic connections on $X$ with rational parabolic weights along $D$.   
\end{abstract}

\baselineskip=15.5pt 

\date{Last updated on \today\,at \currenttime\,(IST)}

\maketitle 

\tableofcontents

\section{Introduction}\label{sec:introduction}
The notion of parabolic vector bundles on a compact Riemann surface $X$ together with a 
finite subset $D \subset X$ of marked points is introduced in \cite{MehSes80} by V. Mehta 
and C.S. Seshadri. In their initial formulation, a parabolic bundle is a vector bundle $E$ 
on $X$ equipped with the data of a flag $F_\bullet(E)$ of $\bb C$--linear subspaces of 
each of its fibers $E_x$ over the points $x \in D$, and a set of increasing real numbers 
in $[0, 1)$ associated to each points of $D$, called weights. The celebrated Mehta--Seshadri 
theorem establishes a one-to-one correspondence between polystable parabolic vector bundles 
on $(X, D)$ of parabolic degree zero and the unitary representations of the fundamental group 
$\pi_1(X\setminus D)$ of the punctured Riemann surface $X\setminus D$. 
Later Simpson reformulated and extended the notion of parabolic vector bundles to the case 
of schemes $X$ with a normal crossings divisor $D$ in \cite{Sim90}. When the weights are 
all rational, and therefore can be assumed to be in $\frac{1}{r}\bb Z$ for some integer $r$, 
this amounts to a locally free coherent sheaf $\mc E$ of finite rank, together with a filtration 
\begin{align*}
	\mc E = \mc E_0\supset \mc E_{\frac{1}{r}}\supset \mc E_{\frac{2}{r}} 
	\supset\cdots\supset \mc E_{\frac{r-1}{r}}\supset \mc E_1= \mc E(-D) 
\end{align*}
by its coherent subsheaves \cite{MarYok92, Bis97}. 
Since their introduction, parabolic bundles, parabolic connections and their 
moduli spaces have been an active area of research in algebraic geometry and 
differential geometry, and have been used in a number of important applications. 
Even though most of the progress has been done in the case of curves, more recently, 
parabolic connections have made appearance in the work of Donagi and Pantev on 
Geometric Langlands conjecture using Simpson's non-abelian Hodge theory \cite{DonPan19}.

When the weights are rational, it is known that parabolic bundles can be 
interpreted equivalently as equivariant bundles on a suitable ramified 
Galois cover \cite{Bis97}. A more intrinsic interpretation is found later using 
vector bundles on certain algebraic stacks associated to $(X, D)$, 
namely the stack of roots in \cite{Bor07}. More precisely, there is a Fourier like 
correspondence between parabolic vector bundles on a scheme and ordinary vector 
bundles on certain stack of roots. This naturally raised the question of 
understanding the parabolic connections on smooth varieties through such 
Fourier like correspondence over root stacks. This has been shown to be 
true over curves \cite{BisMajWon12, LorSaiSim13}, and very recently over 
higher dimensional varieties in \cite{BorLaa23}. This correspondence has been further 
extended to the case of real parabolic connections on a real variety $(X, D)$ in 
\cite{ChPa24}, and for orthogonal and symplectic parabolic connections in \cite{ChMa24}.

In this article we study such correspondences under the presence of a group action. 
More precisely, let $X$ be a connected smooth complex projective variety. 
Let $G$ be a connected complex linear algebraic group acting on $X$. 
Let $D$ be a smooth reduced effective Cartier divisor on $X$ invariant under the 
$G$--action on $X$. Given a positive integer $r$, we can associate a stack of 
$r$-th roots of $\mathcal{O}_X(D)$ and its canonical section vanishing along $D$, 
which we denote by $\mf X$. It comes with a natural morphism 
$\pi : \mf{X} \to X$. The root stack $\mf X$ admits a tautological invertible sheaf  
$\ms{M}$ satisfying $\pi^*\mc{O}_{X}(D) \simeq \ms{M}^{\otimes r}$. 
It is a natural to ask whether the $G$-action on $X$ can be lifted to a $G$-action 
on $\mf X$ making the morphism $\pi : \mc X \to X$ $G$--equivariant. 
We address this question first; we show that, if $G$ has no non-trivial characters, 
then the $G$--action on $(X, D)$ can be lifted to $\mf X$ making $\pi$ $G$-equivariant. 
Moreover, the tautological invertible sheaf  $\ms{M}$ naturally admits a linearization for 
this $G$-action on $\mf X$. This is done in Proposition \ref{prop:lifting-action-on-root-stack} 
and Proposition \ref{prop:linearization-on-tautol-bun-ovr-rt-stack}.

It is shown in \cite{Bor07} that there is an equivalence between the category of algebraic 
vector bundles on $X$ and the category of parabolic vector bundles on $X$ with parabolic 
structures along $D$ and having rational weights in $\frac{1}{r}\mathbb{Z}$. 
We extend this result to the case of $G$-equivariant parabolic bundles 
(see Definition \ref{def:equivariant-parabolic-bundle}).

\begin{theorem}[Theorem \ref{thm:equiv-Biswas-Borne-correspondence}]
	There is an equivalence between the category of $G$-equivariant vector bundles on $\mf{X}$ 
	and the category of $G$-equivariant parabolic bundles on $X$ with parabolic structures 
	along $D$ having rational weights in $\frac{1}{r}\mathbb{Z}$. 
\end{theorem}

Finally, let $\mf D$ be an effective Cartier divisor on $\mf X$ satisfying 
$\mc{O}_{\mf X}(\mf D)\simeq \ms{M}$. The correspondence of N. Borne as above has been 
recently generalized to an equivalence between the category of logarithmic connections 
on $\mf X$ along $\mf D$ and the category of parabolic connections on $X$ along $D$ and 
weights in $\frac{1}{r}$ when $D$ has strict normal crossings \cite{BorLaa23}. 
In the final section, we consider $D$ to be a smooth irreducible effective Cartier divisor 
on $X$, and extend the result of \cite{BorLaa23} to our equivariant setup.

\begin{theorem}[Theorem \ref{thm:equivariant-connection-correspondence}]
	There is an equivalence between the category of $G$-equivariant logarithmic connections on $(\mf{X},\mf{D})$ (Definition \ref{def:equivariant-stacky-logarithmic-connection}) and the category of $G$-equivariant parabolic connections on $(X,D)$ (Definition \ref{def:equivariant-logarithmic-and-parabolic-connection}). 
\end{theorem}

\section{Lifting $G$-action to root stacks}
Unless explicitly specified, all the stacks in this article are algebraic and are defined 
over $\bb C$. We use the same symbol $T$ to denote a $\bb C$-scheme $T$ as well as the 
$\bb C$-stack %$(\Sch\!/T)_{\text{\'et}}$ 
represented by the scheme $T$, if no confusion arises. 
Given $\bb C$-stacks $\mc Y$ and $\mc Y'$, we use the symbol $\mc Y\times\mc Y'$ to 
mean their fiber product $\mc Y\times_{\Spec \bb C}\mc Y'$ over $(\Sch/\bb C)$; note that 
it is a $\bb C$-stack and for any $T \in (\Sch/\bb C)$ we have 
\begin{equation}\label{eqn:fiber-product-valued-points}
	(\mc Y\times\mc Y')(T) = \mc Y(T)\times\mc Y'(T), 
\end{equation}
because $(\Spec \bb C)(T)$ is singleton (c.f. \cite[\S\,4.1]{MartinWise2018} or \cite{Ols07}).

\subsection{Preliminaries on root stack}\label{subsection:root-stack} 
We recall some basic definition and properties of root stack from \cite{Cad07}. 
A \textit{generalized Cartier divisor} on a scheme $T$ is an ordered pair $(M, s)$, 
where $M$ is an invertible sheaf of $\mc O_T$-modules on $T$ and $s \in \Gamma(T, M)$. 
Since a morphism $X \to [\bb A^1/\bb G_m]$ is given by a principal $\bb G_m$-bundle 
on $X$ equipped with a section of the associated $\bb A^1$--bundle on $X$, and vice versa, 
to give a generalized Cartier divisor $(L, s)$ on $X$ is equivalent to give a morphism 
into the quotient stack 
$$f_{(L,\,s)} : X \to [\bb A^1/\bb G_m].$$ 
Given an integer $r \geq 2$, let $\theta_r : [\bb A^1/\bb G_m] \to [\bb A^1/\bb G_m]$ 
be the morphism of quotient stacks induced by the $r$-th power maps on both $\bb A^1$ 
and $\bb G_m$. Then the fiber product 
\begin{equation*}
	{\mf X}_{(L,\, s,\, r)} := X\times_{f_{(L,\, s)},\,[\bb A^1/\bb G_m],\,\theta_r} [\bb A^1/\bb G_m] 
\end{equation*}
is a Deligne-Mumford stack over $\bb C$, called the {\it stack of $r$-th roots of $(L, s)$} 
on $X$, \cite[Theorem 2.3.3]{Cad07}. 
It has the following description: the objects of $\mf X_{(L,\, s,\, r)}$ are given by 
the quadruples $\left(f : U \to X, M, \phi, t\right)$, where 
\begin{itemize}
	\item $f : U \to X$ is a morphism of $\bb C$-schemes, 
	\item $M$ is an invertible sheaf of $\mc O_U$-modules on $U$, 
	\item $t \in \Gamma(U, M)$, and 
	\item $\phi : M^{\otimes r} \stackrel{\simeq}{\longrightarrow} f^*L$ 
	is an $\mc O_U$-module isomorphism such that $\phi(t^{\otimes r}) = f^*s$. 
\end{itemize}
A morphism from 
$\left(f : U \to X, M, \phi, t\right)$ to $\left(f' : U' \to X, M', \phi', t'\right)$ 
in $\mf X_{(L,\, s,\, r)}$ is given by a pair $(g, \psi)$, where 
\begin{itemize}
	\item $g : U \to U'$ is a morphism of $X$-schemes, and 
	\item $\psi : M \to g^* M'$ is an $\mc O_U$-module isomorphism such that $\psi(t)= g^*(t')$. 
\end{itemize}
Then we have a natural morphism of stacks 
\begin{equation}\label{eqn:str-map-of-root-stack-to-X}
	\pi : {\mf X}_{(L,\, s,\, r)} \to X 
\end{equation}
that sends the object $\left(f : U \to X, M, \phi, t\right)$ of $\mf X_{(L,\, s,\,r)}$ 
to the $X$-scheme $f : U \to X$, and a morphism $(g, \psi)$, as above, 
to the morphism of $X$-schemes $g : U \to U'$. 
Note that there is a natural morphism $X \to \Spec\bb C$ (arising from the structure 
morphism of $X$), whose composition with $\pi$ makes $\mf{X}_{(L,\, s,\, r)}$ into 
a $\bb{C}$-stack. 

The projection map onto the second factor 
$$X\times_{f_{(L, s)},\,[\bb A^1/\bb G_m],\,\theta_r}[\bb A^1/\bb G_m] \to [\bb A^1/\bb G_m]$$ 
corresponds to the {\it tautological generalized Cartier divisor} 
%on {\it tautological generalized Cartier divisor} 
$(\ms{M}, \zeta)$ on $\mf X_{(L,\,s,\,r)}$, where $\ms{M}$ is the 
{\it tautological invertible sheaf} on $\mf X$ and $\zeta$ is the 
{\it tautological global section} of $\ms{M}$. 
Note that $\pi^*(L) \cong \ms M^{\otimes r}$.

\subsection{$G$-action on a root stack}\label{sec:G-action-on-a-root-stack}

Let $G$ be a connected affine algebraic group over $\bb C$ acting on a smooth complex projective 
variety $X$ by $\sigma : G\times X \to X$. Let $p_2 : G\times X \to X$ be the projection map 
onto the second factor. Recall that a \textit{linearization} of the $G$-action $\sigma$ on $X$ 
to an invertible sheaf $L$ on $X$ is an isomorphism 
$$\Phi : p_2^*L \stackrel{\simeq}{\longrightarrow} \sigma^*L$$ 
of sheaves of $\mc O_{G\times X}$-modules on $G\times X$ satisfying the following 
cocycle condition coming from compatibility with the $G$-action on $X$: 
\begin{equation*}
	(\mu\times\Id_X)^*\Phi = (\Id_G\times\sigma)^*\Phi\circ p_{23}^*\Phi,
\end{equation*}
where $p_{23} : G\times G\times X \to G\times X$ is the projection morphism onto the 
second and third factors, and $\mu : G\times G \to G$ is the group operation on $G$ 
(see \cite{MumfordGIT}). For $G$ connected, the above cocycle condition 
is automatic for an invertible sheaf $L$ on $X$, see \cite[\S\,7.2, Lemma 7.1]{Dol}. 

Let $D$ be an effective Cartier divisor on $X$ which is invariant under the $G$-action, 
in the sense that $\sigma^*(D) = p_2^*(D)$. Note that, this implies, in particular, that 
$g^*D = D, \ \forall\ g\in G$, under the natural identification $\{g\}\times X \simeq X$. 
Let $s_D \in H^0(X, \mc O_X(D))$ be the canonical section of the invertible sheaf $\mc O_X(D)$ 
on $X$ whose divisor of zeros $Z(s_D)$ equals $D$ (cf. \cite[p.~157]{HartshorneAG}). 
Assume that $G$ has no non-trivial characters. This is the case, for example, when $G = [G, G]$. 
For example, any connected semisimple algebraic group $G$ over $\bb C$ satisfies $G = [G, G]$; 
and standard examples of such groups include 
${\rm SL}_n(\bb{C})$, ${\rm SO}_n(\bb{C})$, ${\rm Sp}_{2n}(\bb{C})$ etc.

\begin{lemma}\label{lem:linearization-on-O_X(D)}
	Suppose that $G$ has no non-trivial characters. Then with the above assumptions, 
	there is a natural $G$--linearization 
	$$\Phi : p_2^*\mc{O}_X(D)\stackrel{\simeq}{\longrightarrow}\sigma^*\mc{O}_X(D)$$ 
	of $\sigma$ on $\mc{O}_X(D)$ such that $\Phi(p_2^*s_D) = \sigma^*s_D$. 
	Consequently, $\mc O_X(nD)$ admits a natural $G$--linearization $\Phi_n$, 
	for all $n \in \bb Z$. 
\end{lemma}

\begin{proof}
	The assumption $p_2^*D = \sigma^*D$ gives rise to a natural isomorphism 
	\begin{align*}
		\Phi : p_2^*\mc{O}_X(D)\stackrel{\simeq}{\longrightarrow} \sigma^*\mc{O}_X(D)
	\end{align*}
	of invertible sheaves on $G\times X$. We need to show that $\Phi(p_2^*s_D) = \sigma^*s_D$. 
	For each $g\in G$, the natural identification $\{g\}\times X\simeq X$ gives rise to an 
	isomorphism of invertible sheaves 
	$$\Phi_g : \mc{O}_X(D) \stackrel{\simeq}{\longrightarrow} \sigma_g^*\mc{O}_X(D)$$ 
	on $X$, where $\sigma_g : X \to X$ is the multiplication by $g$ map. Note that, 
	the sections $\Phi_g(s_D)$ and $\sigma_g^*(s_D)$ of $\sigma_g^*\mc{O}_X(D)$ satisfies 
	\begin{align*}
		Z(\Phi_g(s_D))= Z(s_D) = D \ \ \text{and} \ \ Z(\sigma_g^*(s_D)) = \sigma_g^*(D). 
	\end{align*}
	Since $\sigma_g^*(D) = D$, we see that the sections $\Phi_g(s_D)$ and $\sigma_g^*(s_D)$ 
	have the same divisor of zeros. Since $X$ is a smooth projective variety, there exists 
	a $\lambda_g\in \bb{C}^*$ such that $\Phi_g(s_D)= \lambda_g\cdot \sigma_g^*(s_D)$ 
	(see \cite[Proposition 7.7]{HartshorneAG}). 
	This gives a map $\chi : G \to \bb{C}^*$ defined by $g \mapsto \lambda_g$. 
	It follows from the cocycle condition for the $G$--linearization $\Phi$ on $\mc O_X(D)$ 
	that $\chi$ is a group homomorphism. 
	Since $G$ has no non-trivial characters by assumption, we conclude that $\chi$ is trivial, 
	i.e. $\lambda_g = 1, \ \forall\ g \in G$. From this, the result follows 
	(see \cite[p.~105]{Dol} and \cite[p.~32]{MumfordGIT}). 
\end{proof}

Henceforth we always assume that $G$ has no non-trivial characters, and 
$\mf X := \mf X_{(\mc O_X(D),\, s_D,\, r)}$ be the stack of $r$-th roots of the 
generalized Cartier divisor $(\mc O_X(D), s_D)$ on $X$. 
Let $\pi : \mf X \to X$ be the natural morphism of stacks over $\bb C$ 
as defined in \eqref{eqn:str-map-of-root-stack-to-X}. 

\begin{proposition}\label{prop:lifting-action-on-root-stack}
	With the above assumptions, the $G$-action $\sigma$ on $X$ lifts to a $G$-action 
	$$\sigma_{\mf X} : G\times {\mf X} \rightarrow \mf X$$ 
	on $\mf X$ (c.f. \cite{Romagny-2005}) such that the diagram 
	\begin{equation}\label{diag:lift-of-group-action}
		\begin{gathered}
			\xymatrix{
				G\times\mf{X} \ar[rr]^{\sigma_{\mf X}} \ar[d]_{\Id\times\pi} & & \mf{X} \ar[d]^{\pi}\\
				G\times X \ar[rr]^{\sigma} && X
			}
		\end{gathered}
	\end{equation}
	is $2$-Cartesian. 
\end{proposition}

\begin{proof}
%	Let us denote $\sigma_{\mf{X}}$ by $F$ for simplicity. 
	Thinking of both $G\times \mf {X}$ and $\mf{X}$ as fibered categories over $(\Sch/\bb{C})$, 
	we first describe where $\sigma_{\mf X}$ sends the objects and morphisms of $G\times \mf{X}$\ . 
	For $T \in (\Sch/\bb{C})$, a typical object of $(G\times \mf{X})(T)$ is given by an ordered pair 
	$(g, \tau)$, where $g \in G(T)$ and $\tau = (f : T \to X, M, \varphi, t) \in \mf X(T)$, 
	see \eqref{eqn:fiber-product-valued-points} and \S\,\ref{subsection:root-stack}. 
%	Since $\varphi : M^{\otimes r} \longrightarrow f^*\mc{O}_X(D)$ is an isomorphism of $\mc O_T$-modules 
%	satisfying 
%	\begin{align}\label{eqn:eqn-1}
%		\varphi(t^r) = f^*(s_D), 
%	\end{align}
	Let $h := (g, f) : T \to G\times X$. Then we have 
	\begin{align*}
		p_1\circ h = g \ \ \text{and} \ \ p_2\circ h = f, 
	\end{align*}
	where $p_1 : G\times X \to G$ and $p_2 : G\times X \to X$ are the projection maps onto the 
	first and the second factors, respectively. 
%	Define 
%	\begin{align}\label{eqn:h-sigma}
%		h^{\sigma} := \sigma\circ h : T\rightarrow X\ .
%	\end{align}
	Consider the composite isomorphism of invertible sheaves on $T$ 
%	We have the following chain of isomorphisms of invertible sheaves on $T$: 
	\begin{align*}
		M^{\otimes r} \stackrel{\varphi}{\longrightarrow} 
		f^*\mc{O}_{X}(D) = h^*(p_2^*{\mc O}_X(D)) 
		\stackrel{h^*(\Phi)}{\longrightarrow} 
		h^*(\sigma^*{\mc O}_X(D)) = ({\sigma\circ h})^*\mc{O}_X(D),
	\end{align*}
	where by equality we mean the canonical identification. 
	Note that $(h^*\Phi)\circ\varphi$ sends the section $t^r \in \Gamma(T, M^{\otimes r})$ to 
	\begin{align*}
		\left((h^*\Phi)\circ\varphi\right) (t^r) = (h^*\Phi)(f^*s_D) 
		&= (h^*\Phi)(h^*(p_2^*s_D)) \\ 
		&= h^*(\Phi(p_2^*s_D)) \\
		&= (\sigma\circ h)^*s_D, 
%		&= (h^{\sigma})^*(s_D)\ .
	\end{align*}
	where the last equality follows from Lemma \ref{lem:linearization-on-O_X(D)}. 
	Therefore, the quadruple 
	$$({\sigma\circ h} : T \to X,\, M,\, (h^*\Phi)\circ\varphi,\, t)$$ 
	defines an object in $\mf{X}(T)$. 
	Therefore, we define 
	\begin{align}\label{eqn:action-description-on-objects}
		\sigma_{\mf X}(g, \tau) := ({\sigma\circ h} : T \to X,\, M,\, (h^*\Phi)\circ\varphi,\, t) \in \mf{X}(T). 
	\end{align}
%	Let us now describe the effect of $F$ on arrows. 
	We now describe the effect of $\sigma_{\mf X}$ on the arrows. 
	Let $(g, \tau)$ and $(g', \tau')$ be two objects in $G\times\mf{X}$. %\eqref{eqn:fiber-product-valued-points} 
	Then $g \in G(T)$ and $g' \in G(T')$, for some $\bb C$-schemes $T$ and $T'$, respectively, and that 
	\begin{align}\label{eqn:description-of-tuples}
		\tau := (f : T \to X, M, \varphi, t) \ \text{ and } \ 
		\tau' := (f' : T' \to X, M', \varphi', t'), 
	\end{align} 
	see \S\,\ref{subsection:root-stack}. 
	Then a morphism from $(g, \tau)$ to $(g', \tau')$ in $G\times\mf X$ is given by a triple 
	\begin{align}\label{eqn:arrows-1}
		(\phi,m,(\ell,\theta))
	\end{align}
	where 
	\begin{itemize}
		\item $\phi: T\rightarrow T'$ is a morphism of $\bb C$-schemes, 
		\item $m : T \to T'$ is a morphism of $\bb C$-schemes making the following diagram commutative 
		\begin{equation}\label{diagram:1} 
			\begin{gathered}
				\xymatrix{
					T \ar[rr]^-m \ar[dr]_-{g} && T' \ar[ld]^-{g'} \\
					& G & 
				}
			\end{gathered}
		\end{equation}
	\item $(\ell,\theta) : \tau \rightarrow \tau'$ is a morphism in $\mf{X}$; 
	i.e., $\ell : T \to T'$ is a morphism of $\bb C$-schemes making the 
	following diagram commutative 
	\begin{align}\label{diagram:2}
		\begin{gathered}
		\xymatrix{
			T \ar[rr]^{\ell} \ar[dr]_-{f} & & T'\ar[ld]^-{f'} \\
			& X &
		}
	\end{gathered}
	\end{align}
	and $\theta : M \longrightarrow \ell^*(M')$ is an isomorphism of $\mc O_T$-modules 
	such that $\theta(t) = \ell^*(t')$. 
	\end{itemize}
	Since both $m$ and $(\ell, \theta)$ maps to $\phi$ under the morphism $\pi : \mf X \to X$, 
	it follows that $m = \ell$ in the triple \eqref{eqn:arrows-1} 
	(c.f. \cite[Definition 4.2]{MartinWise2018}). 
	Consider the morphisms 
	\begin{align*}
		h:= (g,f) : T\rightarrow G\times X \ \ \text{and}\ \ h':=(g',f'): T'\rightarrow G\times X\ .
	\end{align*}
	It follows from the commutativity of the diagrams \eqref{diagram:1} and \eqref{diagram:2} and 
	the equality $\ell = m$ that the following diagram commutes. 
	\begin{align*}
		\xymatrix{
			T \ar[rr]^{\ell} \ar[dr]_{h} & & T'\ar[ld]^{h'}\\
			& G\times X &
		}
	\end{align*}
	Therefore, $\sigma\circ h = (\sigma\circ h')\circ\ell$. 
	Thus we get a morphism in $\mf{X}$ given by
	\begin{align*}
		\sigma_{\mf X}(g, \tau) = ({\sigma\circ h}, M, (h^*\Phi)\circ\varphi, t) 
		\xrightarrow{(\ell,\theta)} 
		({\sigma\circ h'}, M', ({h'}^*\Phi)\circ\varphi', t')  
		= \sigma_{\mf X}(g', \tau'),  %\quad (\text{cf.}\ \eqref{diagram:2})
	\end{align*}
	which we define to be $\sigma_{\mf X}((\phi, m, (\ell, \theta)))$. 
	Now it is straight-forward to check that $\sigma_{\mf X}$ is indeed a morphism of stacks 
	making the diagram \eqref{diag:lift-of-group-action} $2$-commutative. 
	
	It remains to show that the diagram \eqref{diag:lift-of-group-action} is $2$-Cartesian. 
	The objects of $\left(G\times X\right)\underset{\sigma,X,\pi}{\times} \mf{X}$ are given by 
	quadruples $(T, h, \tau, \alpha)$, where 
	\begin{itemize}
		\item $T$ is a $\bb{C}$-scheme, 
		\item $h \in (G\times X)(T) = \Hom_{(\Sch/\bb{C})}(T, G\times X)$, 
		\item $\tau \in \mf{X}(T)$, and 
		\item $\alpha : \sigma\circ h \stackrel{\simeq}{\longrightarrow} \pi\circ\tau$ 
		in $X(T) = \Hom_{(\Sch/\bb C)}(T, X)$.
	\end{itemize}
	Since $X(T)$ is a set, it follows that $\alpha$ is identity, i.e., $\sigma\circ h = \pi\circ\tau$. 
	Thus the quadruple $(T, h, \tau, \alpha)$ simplifies to $(T, h, \tau)$ satisfying 
	\begin{align}\label{eqn:eqn-2}
		\sigma\circ h = \pi\circ\tau. 
	\end{align}
	If $h = (g, k)$, for some $g : T \to G$, $k : T \to X$, 
	and if $\tau = (f, M, \varphi, t)$ (see \S\,\ref{subsection:root-stack}), 
%	Note that $\pi\circ\tau = f$. 
	denoting by $\iota : G \to G$ the inversion morphism of $G$, the above relation \eqref{eqn:eqn-2} gives 
	\begin{align*}
		& \sigma\circ(g\ ,k) = f \\
		\implies & \sigma\circ\left(\iota\circ g, \sigma\circ(g, k)\right) = \sigma\circ(\iota\circ g, f), \\
		\implies & k = \sigma\circ(\iota\circ g, f), 
	\end{align*}
	where the last equality follows from the following commutative diagram:
	\begin{align*}
		\xymatrix{T \ar[rr]^(.4){(g,g,k)} \ar[rrdd]_{(e,k)} && 
			G\times( G\times X) \ar[rr]^{\iota\times \sigma} 
			\ar[d]^{\iota\times\Id_G\times\Id_X} && G\times X \ar[dd]^{\sigma}\\
			&& G\times G\times X \ar[d]^{\mu\times \Id_X}&&\\
			&& G\times X \ar[rr]^{\sigma}&& X
		}
	\end{align*}
	where $\mu: G\times G \to G$ is the group multiplication, and 
	$e : T\rightarrow G$ is the constant map sending everything to the neutral point of $G$. 
	Thus we can conclude that the data of $(T, g, \tau)$ determines map $h$, 
	and thus there is an isomorphism of fibered categories 
	\begin{align*}
		\left(\left(G\times X\right)\underset{\sigma, X, \pi}{\times}\mf{X}\right)(T) \simeq & (G\times\mf{X})(T) \\
		\text{sending}\quad (h,\tau)\ \  \mapsto & \ \ (p_1\circ h, \tau)\\
		\text{and}\quad\left((g', \sigma\circ(\iota\circ g', \pi\circ\tau')), \tau'\right)\ \  \mapsfrom & \ \  (g', \tau'),
	\end{align*}
	which is functorial in $T$. Hence the diagram \eqref{diag:lift-of-group-action} is $2$-Cartesian, as required. 
\end{proof}

\begin{remark}\label{rem:second-projection-notation}
	To avoid confusion, we denote by $p_2 : G\times X \to X$ and 
	$q_2 : G\times\mf{X} \to \mf{X}$ the second projections, respectively. 
	Note that, the following diagram is Cartesian. 
	\begin{align*}
		\xymatrix{
			G\times\mf{X} \ar[rr]^-{q_2} \ar[d]_-{\Id\times\pi} && \mf{X} \ar[d]^{\pi} \\
			G\times X \ar[rr]^-{p_2} && X 
		}
	\end{align*} 
\end{remark}

\subsection{$G$--action on the tautological invertible sheaf on the root stack}\label{subsec:G-action-on-tautological-line-bundle-on-root-stack}
%Fix an integer $r \geq 2$, and let $\mf X := \mf{X}_{(\mc O_X(D),\, s_D,\, r)}$ 
%be the stack of $r$-th roots of the generalized Cartier divisor $(\mc O_X(D), s_D)$ 
%on $X$. Note that $\mf X$ is isomorphic to the fiber product 
%$X\times_{t,\,[\bb A^1/\bb G_m],\,\theta_r}[\bb A^1/\bb G_m]$, 
%where $t : X \to [\bb A^1/\bb G_m]$ is the morphism give by the generalized Cartier divisor 
%$(\mc O_X(D), s_D)$ on $X$ and $\theta_r : [\bb A^1/\bb G_m] \to [\bb A^1/\bb G_m]$ is the 
%morphism induced by the $r$-th power maps on both $\bb A^1$ and $\bb G_m$. 
%The projection map $X\times_{t,\,[\bb A^1/\bb G_m],\,\theta_r}[\bb A^1/\bb G_m] \to [\bb A^1/\bb G_m]$ onto the second factor corresponds to a generalized Cartier divisor  
%%on {\it tautological generalized Cartier divisor} 
%$(\ms{M}, \zeta)$ on $\mf X$, where $\ms{M}$ is the tautological invertible sheaf 
%on $\mf X$ and $\zeta$ is the tautological global section of $\ms{M}$; see \cite{Cad07}. 

Let $u : U \to \mf{X}$ be an atlas given by the quadruple 
$(f_u : U \to X, M_u, \varphi_u, t_u)$; see \S\,\ref{subsection:root-stack}. 
The tautological invertible sheaf $\ms{M}$ on $\mf{X}$ is given by the 
invertible sheaf $M_u$ on each such atlas $U$; see \cite[\S\,3.2]{Bor07}. 
We show that the $G$--action on $\mf{X}$, as constructed in 
Proposition \ref{prop:lifting-action-on-root-stack}, can be lifted to a 
$G$--linearization on $\ms{M}$.

\begin{proposition}\label{prop:cartesian-diagram-for-root-stack}
	Let $U$ be a $\bb C$-scheme and let $u : U \to \mf{X}$ be an atlas of $\mf X$ 
	given by the quadruple $(f_u : U \to X, M_u, \varphi_u, t_u)$. 
	Then the Cartesian diagram of schemes 
	\begin{equation*}\label{diag:cartesian-diagram-1}
		\xymatrix{
			U' \ar[rr]^-{\sigma'} \ar[d]_{h_u} && U\ar[d]^{f_u} \\
			G\times X \ar[rr]^-{\sigma} && X 
		}
	\end{equation*}
	induces a Cartesian diagram of the form
	\begin{equation}\label{diag:Cartesian-diagram-for-stack}
		\begin{gathered}
			\xymatrix{
				U' \ar[rr]^{\sigma'} \ar[d]_{u'} && U \ar[d]^{u} \\
				G\times\mf{X} \ar[rr]^{\sigma_{\mf{X}}} && \mf{X}. 
			}
		\end{gathered}
	\end{equation}
\end{proposition}

\begin{proof}
	We first describe the left vertical map $u'$. Suppose $h_u : U' \to G\times X$ 
	is given by the pair of morphisms $(\alpha_u : U' \to G, \beta_u : U' \to X)$. 
	Consider the following composition of isomorphisms 
	\begin{align}\label{eqn:iso-1}
		\begin{gathered}
		\xymatrix{
		(\sigma')^*(M_u^{\otimes r}) \ar[rr]^-{(\sigma')^*\varphi_u}_-{\simeq} &&  (\sigma')^*(f_u^*(\mc{O}_X(D))) \ar[d]_{\simeq} && \beta_u^*(\mc{O}_X(D)) \\ 
		&& h_u^*(\sigma^*\mc{O}_X(D)) \ar[rr]^-{h_u^*{\Phi^{-1}}}_-{\simeq} 
		&& h_u^*(p_2^*\mc{O}_X(D)) \ar[u]_{\simeq} 
	}
	\end{gathered}
	\end{align}
	where the isomorphism at the bottom line of \eqref{eqn:iso-1} 
	follows from Lemma \ref{lem:linearization-on-O_X(D)}. 
	Note that 
	\begin{align}
		& \varphi_u(t^{\otimes r})=f_u^*(s_D) \nonumber \\ 
		\implies & (\sigma')^*(\varphi_u(t^{\otimes r})) = (\sigma')^*(f_u^*(s_D)) \nonumber\\
		\implies & ((\sigma')^*\varphi_u)((\sigma')^*t_u^{\otimes r}) 
		= h_u^*(\sigma^*(s_D)) \nonumber \\ 
		\implies & ((\sigma')^*\varphi_u)((\sigma')^*t_u^{\otimes r}) 
		= h_u^*\left(\Phi(p_2^*(s_D))\right), 
		\quad\text{ by Lemma \ref{lem:linearization-on-O_X(D)}} \nonumber \\
		\implies & ((\sigma')^*\varphi_u)((\sigma')^*t_u^{\otimes r}) 
		= h_u^*(\Phi)(h_u^*(p_2^*(s_D))) \nonumber \\
		\implies & (h_u^*{\Phi^{-1}})\left[((\sigma')^*\varphi_u)((\sigma')^*t_u^{\otimes r})\right] 
		= h_u^*(p_2^*(s_D)) = \beta_u^*(s_D). \label{eqn:section-equality-1}
	\end{align}
	Therefore, by \eqref{eqn:iso-1} and \eqref{eqn:section-equality-1} the quadruple 
	\begin{align*}
		u_0 := \left(\beta_u : U' \to X, (\sigma')^*(M_u), 
		h_u^*(\Phi^{-1})\circ(\sigma')^*\varphi_u, (\sigma')^*(t_u)\right)
	\end{align*}
	defines a morphism $$u_0 : U' \to \mf{X}.$$ 
	
	Now the pair of morphisms 
	$$(\alpha_u : U' \to G,\, u_0 : U' \to \mf{X})$$ 
	defines a morphism $u' : U' \to G\times\mf{X}$. 
	We show that this $u'$ does the job. 
	Note that, the morphism $u\circ\sigma' : U' \to \mf{X}$ is given by the quadruple 
	\begin{align}\label{eqn:tuple-1}
		(f_u\circ\sigma', \left(\sigma')^*(M_u), (\sigma')^*\varphi_u, (\sigma')^*t_u\right), 
	\end{align}
	and the composite morphism $\sigma_{\mf X}\circ u' : U' \to \mf{X}$ is given by 
	the quadruple 
	\begin{equation*}
		\left(\sigma\circ h_u, (\sigma')^*(M_u), 
		(h_u^*\Phi)\circ(h_u^*(\Phi^{-1})\circ(\sigma')^*\varphi_u), (\sigma')^*t_u\right), 
	\end{equation*}
	(c.f. \eqref{eqn:action-description-on-objects}), which is equal to 
	\begin{equation}\label{eqn:tuple-2}
		\left(f_u\circ\sigma', (\sigma')^*(M_u), (\sigma')^*\varphi_u, (\sigma')^*t_u\right)  
	\end{equation}
	since $h_u^*(\Id) = \Id$. 
	Since the quadruples in \eqref{eqn:tuple-1} and \eqref{eqn:tuple-2} are the same, 
	it follows that $$u\circ\sigma' = \sigma_{\mf{X}}\circ u',$$ 
	i.e., the diagram \eqref{diag:Cartesian-diagram-for-stack} commutes. 
	
	To show that the diagram \eqref{diag:Cartesian-diagram-for-stack} is Cartesian, 
	note that the objects of 
	$(G\times \mf{X})\underset{\sigma_{\mf{X}},\mf{X},u}{\times} U$ 
	are quadruples of the form $(T, \gamma, x, \delta)$, where 
	\begin{itemize}
		\item $T$ is a $\bb{C}$-scheme, 
		\item $\gamma \in (G\times \mf{X})(T)$, 
		\item $x \in U(T) = \Hom_{(\Sch/\bb C)}(T, U)$, and 
		\item $\delta : \sigma_{\mf{X}}(\gamma) \stackrel{\simeq}{\longrightarrow} u(x)$ 
		is an isomorphism in $\mf{X}(T)$. 
	\end{itemize}
	
	Let $\gamma = (\alpha, \beta) \in (G\times\mf X)(T)$, for some $\alpha: T \to G$ 
	and $\beta: T \to \mf{X}$. Suppose that $\beta$ is given by the quadruple 
	$(f_\beta: T \to X, M_{\beta}, \varphi_{\beta}, t_{\beta})$, and let 
	$h_{\beta} := (\alpha, f_{\beta}) : T \to G\times X$. 
	Then we have 
	\begin{align}\label{eqn:tuple-3}
		\sigma_{\mf{X}}(\gamma) = \left(\sigma\circ h_{\beta}, M_{\beta}, h_{\beta}^*(\Phi)\circ\varphi_{\beta}, t_{\beta}\right); 
	\end{align}
	see \eqref{eqn:action-description-on-objects}. 
	Since $u : U \to \mf X$ is given by the quadruple $(f_u : U \to X, M_u, \varphi_u, t_u)$, 
	it follows that $u(x) := u\circ x$ is given by the quadruple 
	\begin{align}\label{eqn:tuple-4}
		\left(f_u\circ x, x^*(M_u), x^*(\varphi_u), x^*(t_u)\right). 
	\end{align}
	Since any morphism in the groupoid $\mf{X}(T)$ maps to the identity map $\Id_T$ under $\pi$, 
	it follows from \eqref{eqn:tuple-1} and \eqref{eqn:tuple-2} that 
	$\sigma\circ h_{\beta} = f_u\circ x$, i.e., the following diagram commutes:
	\begin{align*}
		\xymatrix{
			T\ar[rr]^-{x} \ar[d]_{h_\beta} && U \ar[d]^{f_u} \\
			G\times X \ar[rr]^-{\sigma} && X}
	\end{align*}
	Thus there exists a unique morphism 
	$T\rightarrow U' = (G\times X)\underset{\sigma,\,X,\,u}{\times} U$. 
	This gives a bijection
	\begin{align*}
		\left((G\times \mf{X})\underset{\sigma_{\mf{X}},\,\mf{X},\,u}{\times}U\right)(T)\simeq U'(T)
	\end{align*}
	which is functorial in $T$, proving that the 
	diagram \eqref{diag:Cartesian-diagram-for-stack} is Cartesian.
\end{proof}

\begin{corollary}\label{cor:pullback-of-an-atlas-is-an-atlas}
	For any atlas $u : U \to \mf{X}$ of $\mf X$ with $U$ a $\bb C$-scheme, 
	its pullback $u' : U' \to G\times \mf{X}$, as defined in 
	Proposition \ref{prop:cartesian-diagram-for-root-stack}, 
	is an atlas of $G\times\mf X$. 
\end{corollary}

\begin{proof}
	To show that $u'$ gives an atlas, we need to check that for any scheme 
	$t : T \to G\times\mf{X}$, the projection map 
	$T \underset{t,\,G\times\mf{X},\,u'}{\times}U' \to T$ 
	is smooth and surjective \cite[Definition 8.1.4]{Ols16}. 
	It follows from Proposition \ref{prop:cartesian-diagram-for-root-stack} that 
	\begin{align*}
		T\underset{t,\,G\times\mf{X},\,u'}{\times} U' 
		& \simeq T \underset{t,\,G\times\mf{X},\,u'}{\times} 
		\left((G\times\mf{X})\underset{\sigma_{\mf{X}},\,\mf{X},\,u}{\times} U\right) \\ 
		& \simeq T\underset{(\sigma_{\mf{X}}\circ t),\,\mf{X},\,u}{\times} U. 
	\end{align*}
	Since $U\overset{u}{\longrightarrow} \mf{X}$ is an atlas, the morphism 
	$T\underset{(\sigma_{\mf{X}}\circ t),\,\mf{X},\,u}{\times} U \rightarrow T$ 
	is a smooth surjection. Hence the result follows. 
\end{proof}

\begin{lemma}\label{lem:action-morphism-on-root-stack-is-smooth}
	With the above notations, $\sigma_{\mf{X}} : G\times\mf X \to \mf X$ 
	is a smooth morphism of algebraic stacks. 
\end{lemma}

\begin{proof}
	We need to check that for some atlases $v : V \to G\times \mf{X}$ 
	and $u : U \to \mf{X}$, there exists a commutative diagram
	\begin{align*}
		\xymatrix{
			V \ar[rr]^-{f} \ar[d]_{v} && U \ar[d]^{u} \\
			G\times\mf{X} \ar[rr]^-{\sigma_{\mf{X}}} && \mf{X}
		}
	\end{align*}  
	such that $f$ is smooth (see \cite[Definition 2.4]{Heinloth2010}). 
	Choose an atlas $u : U \to \mf{X}$ given by the quadruple 
	$(f_u, M_u, \varphi_u, t_u)$, and consider the following Cartesian diagram 
	\begin{align*}
		\xymatrix{
			U' \ar[rr]^-{\sigma'} \ar[d]_{h_u} && U \ar[d]^{f_u} \\ 
			G\times X \ar[rr]^-{\sigma} && X.
		}
	\end{align*}
	Since the action morphism $\sigma : G\times X \to X$ is smooth, 
	so is its pullback $\sigma' : U' \to U$. 
	%\cite[Proposition 10.1]{HartshorneAG}. 
	It follows from Corollary \ref{cor:pullback-of-an-atlas-is-an-atlas} 
	that the pullback $u' : U' \to G\times \mf{X}$ of $u$ is an 
	atlas on $G\times\mf{X}$. Thus, we have a commutative diagram
	\begin{align*}
		\xymatrix{U' \ar[rr]^-{\sigma'} \ar[d]_{u'} && U \ar[d]^{u}\\
			G\times \mf{X} \ar[rr]^-{\sigma_{\mf{X}}} && \mf{X}
		}
	\end{align*} 
	where $u$ and $u'$ are atlases and $\sigma'$ is smooth. 
	Thus $\sigma_{\mf{X}}$ is a smooth morphism of stacks.
\end{proof}

Since $\sigma_{\mf{X}} : G\times\mf X \to \mf X$ is a smooth morphism, 
it induces a morphism of Lisse-\'etale sites 
\begin{align*}
	\text{Lis-\'Et}(G\times \mf{X}) \longrightarrow \text{Lis-\'Et} (\mf{X})
\end{align*}
(see \cite[\href{https://stacks.math.columbia.edu/tag/0GR1}{Tag 0GR1}]{StackProject} 
and \cite{Ols07}), 
which enables us to define pull-back $\sigma_{\mf{X}}^*(\mc{F})$ of a 
coherent sheaf $\mc{F}\in \mf{X}_{\text{lis-\'et}}$ as follows. 
Given a smooth morphism $t : T \to G\times\mf{X}$, the composition  
$\sigma_{\mf{X}}\circ t : T \to \mf{X}$ is smooth by 
Lemma \ref{lem:action-morphism-on-root-stack-is-smooth}. 
Following notation and proof of \cite[Lemma 6.5]{Ols07}, we define 
\begin{align*}
	(\sigma_{\mf X}^*\mc{F})_t := \mc{F}_{\sigma_{\mf X}\circ t}\,. 
\end{align*}

\begin{proposition}\label{prop:linearization-on-tautol-bun-ovr-rt-stack}
	The tautological invertible sheaf  $\ms{M}$ on $\mf{X}$ admits a linearization of the 
	$G$-action on $\mf{X}$. Consequently, every tensor power of $\ms{M}$ on $\mf X$ 
	also admits a linearization.
\end{proposition}

\begin{proof}
	Let $V$ be a $\bb C$-scheme and let $v : V \to G\times\mf{X}$ be a smooth morphism 
	given by a pair of morphisms $\alpha_v : V \to G$ and $\beta_v : V \to \mf{X}$, 
	where $\beta_v$ is given by the quadruple 
	\begin{align}\label{eqn:tuple-5}
		(f_v : V \to X,\, M_v,\, \varphi_v,\, t_v). 
	\end{align}
	Let $h_v : = (\alpha_v,\, f_v) : V \to G\times X$. 
	Then $\sigma_{\mf{X}}\circ v : V \to \mf{X}$ is given by the quadruple 
	$$\sigma_{\mf X}(\alpha_v, \beta_v) = \left(\sigma\circ(\alpha_v, f_v) : 
	V \to X,\, M_v,\, \big((\alpha_v, f_v)^*\Phi\big)\circ\varphi_v, t_v\right),$$
	(see \eqref{eqn:action-description-on-objects}), which implies that 
	\begin{align*}
		(\sigma_{\mf{X}}^*\ms{M})_{v} = \ms{M}_{\sigma_{\mf{X}}\circ v} = M_v.
	\end{align*}
	Similarly, for the second projection $q_2 : G\times\mf{X} \to \mf{X}$, 
	the map $q_2\circ v  = \beta_v : V \to \mf{X}$ is given by the quadruple 
	\eqref{eqn:tuple-5}, and thus 
	\begin{align*}
		(q_2^*\ms{M})_{v} = \ms{M}_{q_2\circ v} = \ms{M}_{\beta_v} = M_v. 
	\end{align*}
	Thus there exists a natural isomorphism 
	\begin{align*}
		\sigma_{\mf{X}}^*\ms{M}\simeq q_2^*\ms{M}. 
	\end{align*}
	The second statement about tensor powers is obvious. 
\end{proof}

%%%%%%%%%%%%%%%%%%%%%%%
\section{Biswas--Borne correspondence for equivariant bundles}
In \cite{Bis97} it is shown that the category of parabolic bundles on $X$ with parabolic 
structure defined along a simple normal crossing divisor $D$ with rational parabolic weights 
is equivalent to the category of orbifold vector bundles on certain finite Galois cover of $X$. 
The notion of parabolic bundles is reformulated in the language of functors in \cite{Yok95}, 
and the main result of \cite{Bis97} is reformulated and generalized using the language of 
root stacks by N. Borne in \cite{Bor07}. 
We begin with the definition of parabolic bundles from \cite{Yok95}. 

\subsection{Parabolic bundles}
Let $\bb R$ be the category whose objects are real numbers and given $r, s \in \bb R$ 
we define 
$${\rm Mor}_{\bb R}(r, s) = \left\{
\begin{array}{ccl}
	\{\iota^{r, s}\}, & \text{ if } & r \leq s, \\ 
	\emptyset, & \text{ if } & r > s. 
\end{array}
\right.$$
Let $\bb R^{\rm op}$ be the opposite category of $\bb R$. 
An {\it $\bb R$-filtered $\mc O_X$-module} is a functor 
$E_\bullet : {\bb R}^{\rm op} \to \Mod(\mc O_X)$, 
where $\Mod(\mc O_X)$ is the category of $\mc O_X$-modules on $X$. 
We denote by $E_t$ the $\mc O_X$-module $E_\bullet(t)$, for all $t \in \bb R$. 
For $s \geq t$ in $\bb R$, we denote by $\iota_{E_\bullet}^{s, t}$ the $\mc O_X$-module 
homomorphism 
$$\iota_{E_\bullet}^{s, t} := E_\bullet(\iota^{s, t}) : E_s \to E_t.$$ 
Morphisms of $\bb R$-filtered $\mc O_X$-modules are given by natural transformations 
of functors. Let $\Filt_{\bb R}(\Mod(\mc O_X))$ be the category whose objects are 
$\bb R$--filtered $\mc O_X$-modules on $X$ and morphisms are as defined above. 
Given $E_\bullet \in \Filt_{\bb R}(\Mod(\mc O_X))$ and $s \in \bb R$, 
we define $$E_\bullet[s] : {\bb R}^{\rm op} \to \Mod(\mc O_X)$$ by setting 
\begin{equation*}
	E_\bullet[s](t) := E_\bullet(s+t), \ \ \text{ and } \ \ 
	E_\bullet[s](\iota^{t, t'}) = \iota_{E_\bullet}^{s+t, s+t'},
\end{equation*}
for all $t, t' \in \bb R$. 
For $s \geq t$ in $\bb R$, we have a morphism of functors 
$$\iota_{E_\bullet}^{[s, t]} : E_\bullet[s] \to E_\bullet[t]$$
given by the unique morphism 
$$\iota_{E_\bullet}^{s+u, t+u} : E_{s+u} \to E_{t+u}, \ \forall\ u \in \bb R.$$
Given $E_\bullet \in \Filt_{\bb R}(\Mod(\mc O_X))$ and an $\mc O_X$-module $V$ on $X$, 
we have an $\bb R$-filtered $\mc O_X$-module $E_\bullet\otimes F$ defined by 
$(E_\bullet\otimes F)(t) := E_t\otimes F, \ \forall\ t \in \bb R$. 
Similarly, tensoring $E_\bullet$ with an $\mc O_X$-module homomorphism 
$f : F \to G$, we get we get a morphism 
$\Id_{E_\bullet}\otimes f : E_\bullet\otimes F \to E_\bullet\otimes G$ 
in $\Filt_{\bb R}(\Mod(\mc O_X))$. 

Fix an integer $r \geq 1$, and let $\frac{1}{r}\bb Z$ be the full subcategory 
of $\bb R$ whose objects are elements of the set $\frac{1}{r}\bb Z$. 
We denote by $\Vect(X)$ the category of {\it vector bundles} 
(locally free coherent sheaves of finite ranks) on $X$. 
Following above notations, the notion of 
{\it $\frac{1}{r}\bb Z$--filtered vector bundles} 
on $X$ makes sense. 

\begin{definition}\label{def:parabolic-bundle}
	A {\it parabolic vector bundle} on $X$ {with a parabolic structure along $D$} 
	with rational parabolic weights in $\frac{1}{r}\bb Z$ is a pair 
	$(E_\bullet, \mf j_{E_\bullet})$, where 
	\begin{equation*}
		E_{\bullet} : \left(\frac{1}{r}\bb{Z}\right)^{\rm op} \longrightarrow {\Vect}(X) 
	\end{equation*}
	is an $\frac{1}{r}\bb Z$--filtered vector bundle on $X$ and 
	\begin{equation*}
		\mf j_{E_\bullet} : E_\bullet\otimes\mc O_X(-D) \longrightarrow E_\bullet[1] 
	\end{equation*}
	is an isomorphism of functors such that the following diagram commutes. 
	\begin{equation}\label{diag:parabolic-bundle-defn} 
		\begin{gathered}
			\xymatrix{
				E_\bullet\otimes\mc O_X(-D) \ar[rr]^-{\mf j_{E_\bullet}}_-{\simeq} 
				\ar[dr]_-{\Id_{E_\bullet}\otimes\iota_D} && E_\bullet[1] 
				\ar[dl]^-{\iota_{E_\bullet}^{[1, 0]}} \\ 
				& E_\bullet & 
			}
		\end{gathered}
	\end{equation}
	where $\iota_D : \mc O_X(-D) \hookrightarrow \mc O_X$ is the natural inclusion map 
	of $\mc O_X$-modules. 
%	
%	In other words, for each integer $n \geq 1$ there is a natural isomorphism of functors 
%	$$E_{\bullet}\otimes\mc{O}_X(-nD) \stackrel{\simeq}{\longrightarrow} E_{\bullet+n}$$
%	such that the following diagram commutes for each $n\geq 0$ and $l\in {\bb Z}$:
%	\begin{equation}\label{diag:parabolic-bundle-defn} 
%		\begin{gathered}
%			\xymatrix{ E_{\frac{l}{r}+n} \ar[r] & E_{\frac{l}{r}} \\
%				E_{\frac{l}{r}}\otimes\mathcal{O}_X(-nD) \ar[u]^{\simeq} \ar[ru] & 
%			}
%		\end{gathered}
%	\end{equation}
	We simply call $E_\bullet$ a parabolic bundle, when the data of $D$ and $r$ 
	is clear from the context. A morphism of parabolic bundles from $E_{\bullet}$ to 
	$F_{\bullet}$ is a natural transformation of functors compatible with the diagram 
	\eqref{diag:parabolic-bundle-defn} in the obvious sense. 
	We denote by ${\PVect}(X, r, D)$ the category of parabolic vector bundles on $X$ having 
	parabolic structure along $D$ and rational parabolic weights in $\frac{1}{r}\bb Z$. 
\end{definition}

Then we have the following result due to N. Borne. 

\begin{theorem}\cite[Th\'eor\`eme 3.13]{Bor07}\label{thm:Biswas-Borne-correspondence} 
	With the above notations, there is an equivalence of categories between $\Vect(\mf{X})$ 
	and ${\PVect}(X, r, D)$ given by the following functors: 
	\begin{itemize}
		\item A vector bundle $\mc E$ on $\mf X$ gives rise to the parabolic vector bundle 
		$E_{\bullet}$ on $X$ defined by the functor 
		$$\frac{\ell}{r} \longmapsto E_{\frac{\ell}{r}} 
		:= \pi_*\left(\mathcal{E}\otimes_{{\mc O}_{\mf X}} 
		\ms M^{-\ell}\right),\,\,\forall\,\ell\,\in \bb Z.$$ 
		
		\item Conversely, a parabolic vector bundle $E_{\bullet}$ on $X$ 
		gives rise to the vector bundle 
		$$\mathcal{E} := \int^{\frac{1}{r}\bb{Z}} 
		\pi^*(E_{\bullet})\otimes_{{\mc O}_{\mf X}} \ms{M}^{\bullet r}$$
		on the root stack $\mf{X}$, where $\ms M$ is the tautological invertible sheaf 
		on ${\mf X}$ satisfying $\pi^*\mc O_X(D) \cong \ms M^{r}$, 
		and $\int^{\frac{1}{r}\mathbb{Z}}$ stands for the coend.
	\end{itemize}
\end{theorem}

\subsection{$G$--equivariance structure}
Let $G$ be an affine algebraic group over $\mathbb{C}$. We assume that $G$ has no 
non-trivial characters (c.f. Proposition \ref{prop:lifting-action-on-root-stack}). 
In this subsection we define the notion of $G$--equivariant parabolic vector bundles on $X$, 
and the notion of $G$--equivariant vector bundles on the root stack $\mf X$, 
and then formulate Biswas--Borne correspondence (Theorem \ref{thm:Biswas-Borne-correspondence}) 
in $G$--equivariant setup. Let 
$$E_{\bullet}: \left(\frac{1}{r}\bb{Z}\right)^{\rm{op}} \rightarrow\Vect(X)$$
be a parabolic vector bundle on $X$. 
If $f : Y \to X$ is a smooth morphism of $\bb C$-schemes, then $f^*D$ is a simple normal 
crossing divisor on $Y$ 
(see \cite[\href{https://stacks.math.columbia.edu/tag/0CBN}{Section 0CBN}]{StackProject}) 
and the functor 
\begin{equation}
	f^*(E_\bullet) : \left(\frac{1}{r}\bb Z\right)^{\rm op} \longrightarrow \Vect(X), 
	\ \ \frac{\ell}{r} \mapsto f^*(E_{\frac{\ell}{r}}), 
\end{equation}
defines a parabolic vector bundle on $Y$ with rational parabolic weights in $\frac{1}{r}\bb Z$ 
and having parabolic structure along $f^*D$, called the {\it pullback of $E_\bullet$ along $f$}. 
Since the $G$--action map $\sigma$ and the second projection map $p_2$ from 
$G\times X$ to $X$ are smooth and $p_2^*(D) = \sigma^*D$ by assumption, both 
$\sigma^*(E_\bullet)$ and $p_2^*(E_\bullet)$ are parabolic bundles on $G\times X$ 
with rational parabolic weights in $\frac{1}{r}\bb Z$ and having parabolic structures along $\sigma^*D = p_2^*D$. 

%{\color{red} {\bf  We should explain how pullback of parabolic bundles and connections are defined in this context. It would be nice if we can prove some result similar to \cite{Alfaya-Biswas-2023} in higher dimensions and in $G$-equivariant setup in a new section. 
%We should write a separate section on parabolic Higgs bundle case and following \cite{Simpson-92} 
%we may show for $G = \bb C^*$ how $\bb C^*$-equivariant parabolic Higgs bundles gives 
%Hodge decomposition and the related picture for the Higgs bundles on root stack.  
%}}

\begin{definition}\label{def:equivariant-parabolic-bundle}
	With the above notations, a {\it $G$-equivariance structure} on a parabolic 
	vector bundle $E_{\bullet}$ is a natural isomorphism of functors 
	\begin{equation*}
		\psi_{\bullet} : p_2^*(E_{\bullet}) \longrightarrow \sigma^*(E_{\bullet}), 
	\end{equation*}
	such that the following diagram commutes: 
	\begin{equation*}
		\begin{gathered}
			\xymatrix{ 
			p_2^*(E_\bullet\otimes\mc O_X(-D)) \ar[d]_-{p_2^*({\mf j}_{E_\bullet})} 
			\ar[rr]^-{\psi_\bullet\otimes\Phi_{-1}} && \sigma^*(E_\bullet\otimes\mc O_X(-D)) 
			\ar[d]^-{\sigma^*({\mf j}_{E_\bullet[1]})} \\ 
			p_2^*(E_\bullet[1]) \ar[rr]^-{\psi_\bullet[1]} && \sigma^*(E_\bullet[1]) 
		}
		\end{gathered}
	\end{equation*}
	where $\Phi_{-1}$ is the isomorphism induced from Lemma \ref{lem:linearization-on-O_X(D)}. 
	A {\it $G$-equivariant parabolic vector bundle} is a pair $(E_\bullet, \psi_\bullet)$, 
	where $E_\bullet$ is a parabolic vector bundle and $\psi_\bullet$ is a $G$-equivariance 
	structure on $E_\bullet$. 
	
	Let $(E_\bullet, \phi_\bullet)$ and $(F_\bullet, \psi_\bullet)$ be two $G$-equivariant 
	parabolic vector bundles on $X$. A morphism from $(E_\bullet, \phi_\bullet)$ to 
	$(F_\bullet, \psi_\bullet)$ is a morphism of parabolic bundles 
	$$f_\bullet : E_\bullet \to F_\bullet$$ 
	such that the following diagram is $2$-commutative. 
	\begin{equation*}
		\xymatrix{
			p_2^*(E_\bullet) \ar[rr]^{p_2^*(f_\bullet)} \ar[d]_-{\phi_\bullet} 
			&& p_2^*(F_\bullet) \ar[d]^-{\psi_\bullet} \\ 
			\sigma^*(E_\bullet) \ar[rr]^-{\sigma^*(f_\bullet)} && \sigma^*(F_\bullet)
		}
	\end{equation*}
	We denote by ${\PVect}^G(X, r, D)$ the category of $G$--equivariant parabolic vector 
	bundles on $X$ with rational parabolic weights in $\frac{1}{r}\bb{Z}$ supported on $D$. 
\end{definition}

Let $\Coh(\mf X)$ be the category of coherent sheaves of $\mc O_{\mf X}$--modules on $\mf X$. 
Let $\sigma_{\mf X} : G\times\mf X \to \mf X$ be the $G$--action on $\mf X$ 
(c.f. Proposition \ref{prop:lifting-action-on-root-stack}), and let 
$q_2 : G\times\mf X \to \mf X$ be the projection morphism onto the second factor.

\begin{definition}\label{def:G-equivariant-bundles-on-G-stack}
	\begin{enumerate}[(a)]
		\item A {\it $G$--linearization} on a coherent sheaf $\mc E$ on $\mf X$ is an isomorphism 
		$$\psi : q_2^*\mc E \stackrel{\simeq}{\longrightarrow} \sigma_{\mf X}^*\mc E$$ 
		of sheaves of $\mc O_{G\times\mf X}$--modules such that the following triangle is 
		$2$--commutative:  
		\begin{equation*}
			\begin{gathered}
				\xymatrix{
					q_{23}^*q_2^*\mc E \ar@{=}[d] \ar[rr]^{q_{23}^*\psi} 
					&& q_{23}^*\sigma^*\mc E \ar@{=}[d] \\ 
					(\mu\times\Id_{\mf X})^*p_2^*\mc E \ar[d]_-{(\mu\times\Id_{\mf X})^*\psi} 
					&& (\Id_G\times\sigma_{\mf X})^*q_2^*\mc E 
					\ar[d]^-{(\Id_G\times\sigma_{\mf X})^*\psi} \\ 
					(\mu\times\Id_{\mf X})^*\sigma_{\mf X}^*\mc E \ar@{=}[rr] 
					&& (\Id_G\times\sigma_{\mf X})^*\sigma_{\mf X}^*\mc E\,, 
				}
			\end{gathered}
		\end{equation*}
		where $\mu : G\times G \to G$ is the group operation on $G$ and 
		$q_{23} : G\times G\times\mf X \to G\times\mf X$ is the 
		projection morphism onto the second and third factors. 
		
		\item A {\it $G$-equivariant coherent sheaf} on $\mf X$ is a pair $(\mc E, \psi)$, 
		where $\mc E$ is a coherent sheaf on $\mf X$ and $\Psi$ is a $G$--linearization 
		on $\mc E$. 
		
		\item A {\it morphism} $(\mc E, \psi) \to (\mc E', \psi')$ of $G$--equivariant 
		coherent sheaves on $\mf X$ is an ${\mc O}_{\mf X}$--module homomorphism 
		$f : \mc E \to \mc E'$ such that the following diagram is $2$--commutative. 
		\begin{equation*}
			\begin{gathered}
				\xymatrix{
					q_2^*\mc E \ar[rr]^-{\psi} \ar[d]_-{q_2^*(f)} && 
					\sigma_{\mf X}^*\mc E \ar[d]^-{\sigma_{\mf X}^*(f)} \\ 
					q_2^*\mc E' \ar[rr]^-{\psi'} && \sigma_{\mf X}^*\mc E' 
				}
			\end{gathered}
		\end{equation*}
	\end{enumerate}
\end{definition}

Let $\Coh^G(\mf{X})$ be the category of $G$--equivariant coherent sheaves on $\mf X$, 
and let $\Vect^G(\mf X)$ be the full subcategory of $\Coh^G(\mf X)$ whose objects are 
$G$--equivariant vector bundle on $\mf{X}$. 

\begin{theorem}\label{thm:equiv-Biswas-Borne-correspondence}
	With the above assumptions, there is an equivalence of categories between 
	${\Vect}^G(\mf{X})$ and ${\PVect}^{G}(X, r, D)$. 
\end{theorem}

\begin{proof} 
	It suffices to show that the equivalence in Theorem \ref{thm:Biswas-Borne-correspondence} 
	preserves the $G$-equivariance structures. 
	Let $(E_\bullet, \psi_\bullet) \in \PVect^G(X, r, D)$ be a $G$--equivariant parabolic 
	vector bundle on $X$, and let 
	$${\mc E} := \int^{\frac{1}{r}\bb Z}\pi^*(E_\bullet)\otimes\ms{M}^{\bullet\,r}$$ 
	be the associated vector bundle on the root stack $\mf X$ obtained from $E_\bullet$ 
	via the coend (certain colimit) construction in the category $\Vect(\mf X)$ of 
	vector bundles on $\mf X$ \cite{Bor07}. 
	Since $q_2^*$ is left adjoint to ${q_2}_*$ (see \cite[Proposition 9.3.6]{Ols16}), 
	it follows that the vector bundle $q_2^*\mc E \in \Vect(G\times\mf X)$ can be interpreted 
	as the following coend: 
	\begin{equation*}
		q_2^*\mc E = \int^{\frac{1}{r}\bb Z} 
		 q_2^*\left(\pi^*(E_\bullet)\otimes\ms{M}^{\bullet\,r}\right). 
	\end{equation*}
	Similarly the vector bundle $\sigma_{\mf X}^*\mc E \in \Vect(G\times\mf X)$ can be 
	interpreted as the following coend: 
	\begin{equation*}
		\sigma_{\mf X}^*\mc E = \int^{\frac{1}{r}\bb Z} 
		\sigma_{\mf X}^*\left(\pi^*(E_\bullet)\otimes\ms{M}^{\bullet\,r}\right)\,. 
	\end{equation*}
	Now, consider the natural isomorphisms of bifunctors 
	$$\left(\frac{1}{r}\bb{Z}\right)^{\rm op}\times\frac{1}{r}\bb{Z} \longrightarrow \Vect(\mf X)$$ 
	given by the following sequence of isomorphisms: 
	\begin{align*}
		q_2^*\left(\pi^*(E_{\bullet})\otimes\ms{M}^{\bullet\,r}\right) & \simeq 
		q_2^*(\pi^*(E_{\bullet}))\otimes q_2^*\ms{M}^{\bullet\,r} \\
		& \simeq (\Id\times\pi)^*(p_2^*(E_{\bullet}))\otimes p_2^*(\ms{M}^{\bullet\,r}) \quad[\text{see Remark}\ \ref{rem:second-projection-notation}]\\
		& \simeq (\Id\times\pi)^*(\sigma_{\mf{X}}^*(E_{\bullet}))\otimes \sigma_{\mf{X}}^*(\ms{M}^{\bullet\ r})\quad[\text{Proposition  \ref{prop:linearization-on-tautol-bun-ovr-rt-stack}}] \\ 
		& \simeq \sigma_{\mf X}^*(\pi^*(E_{\bullet}))\otimes
		\sigma_{\mf X}^*(\ms{M}^{\bullet\,r}) \qquad\quad [\text{from diagram 
			\eqref{diag:lift-of-group-action}}] \\
		& \simeq \sigma_{\mf X}^*(\pi^*(E_{\bullet})\otimes\ms{M}^{\bullet\,r}). 
	\end{align*}
	Then it follows from Remark \ref{rem:isomorphism-of-colimits} (see below) 
	that the coend of $q_2^*\left(\pi^*(E_{\bullet})\otimes\ms{M}^{\bullet\,r}\right)$ 
	and the coend of $\sigma_{\mf X}^*\left(\pi^*(E_{\bullet})\otimes\ms{M}^{\bullet\,r}\right)$ 
	are naturally isomorphic. Thus, we have a natural isomorphism 
	\begin{align}\label{eqn:G-bundle-isomorphism-on-stack}
		\varphi : q_2^*(\mc E) \stackrel{\simeq}{\longrightarrow} \sigma_{\mf X}^*(\mc{E}) 
	\end{align}
	of vector bundles over $G\times\mf X$ which gives the required $G$--equivariance 
	structure on $\mc E$ (see Definition \ref{def:G-equivariant-bundles-on-G-stack}). 
	
	Conversely, let $\mc F \in {\Vect}^G(\mf X)$ be a $G$--equivaraint vector bundle on $\mf X$. 
	Consider the following Cartesian diagram \eqref{diag:lift-of-group-action} from 
	Proposition \ref{prop:lifting-action-on-root-stack}, namely 
	\begin{align*}
		\xymatrix{
			G\times\mf{X} \ar[rr]^-{\sigma_{\mf X}} \ar[d]_-{\Id_G\times\pi} 
			&& \mf{X} \ar[d]^{\pi} \\
			G\times X \ar[rr]^-{\sigma} && X
		}
	\end{align*}
	Since $\pi$ is a finite morphism \cite[Corrolaire 3.6]{Bor07} and $\sigma$ is smooth 
	(in particular flat), we have an isomorphism of functors 
	$\sigma^*\circ\pi_* \simeq (\Id\times\pi)_*\circ\sigma_{\mf{X}}^*$, 
	and similarly, $p_2^*\circ\pi_* \simeq (\Id\times\pi)_*\circ q_2^*$ 
	(see Remark \ref{rem:second-projection-notation}). 
	It follows that for all $\ell\in\bb{Z}$, 
	\begin{align*}
		p_2^*(\pi_*(\mc{F}\otimes_{{\mc O}_{\mf X}}\ms{M}^{-\ell})) 
		& \simeq (\Id\times\pi)_*\left(q_2^*(\mc{F}\otimes_{{\mc O}_{\mf X}}\ms{M}^{-\ell})\right) \\
		& \simeq (\Id\times\pi)_*\left(\sigma_{\mf X}^*(\mc{F}\otimes_{{\mc O}_{\mf X}} 
		\ms{M}^{-\ell})\right), 
		\quad[\text{see Proposition \ref{prop:linearization-on-tautol-bun-ovr-rt-stack}}] \\
		& \simeq \sigma^*\left(\pi_*(\mc{F}\otimes_{{\mc O}_{\mf X}}\ms{M}^{-\ell})\right). 
	\end{align*}
	In other words, if $F_{\bullet}$ is the parabolic bundle corresponding to $\mc{F}$ 
	in Theorem \ref{thm:Biswas-Borne-correspondence}, we obtain isomorphisms
	\begin{align*}
		\phi_{_{\frac{\ell}{r}}} : p_2^*(F_{\frac{\ell}{r}}) 
		\stackrel{\simeq}{\longrightarrow} \sigma^*(F_{\frac{\ell}{r}}), 
		\ \ \forall\ \ell\in\bb{Z}, 
	\end{align*}
	which is functorial in $\ell \in \bb Z$. It is easy to see that the diagrams in 
	Definition \ref{def:equivariant-parabolic-bundle} under these $\phi_{_{\frac{\ell}{r}}}$'s 
	also commute. Thus $F_{\bullet}$ has an induced $G$--equivariance structure.  
\end{proof}

\begin{remark}\label{rem:isomorphism-of-colimits}
	Let $\ms C$ be any category, and let $F, G : J \to \ms C$ be two functors from an 
	index category $J$ such that both the colimits $\lim\limits_{j\in J} F_j$ and 
	$\lim\limits_{j\in J} G_j$ exist in $\mc C$. Then a natural isomorphism 
	of functors $\varphi := \{\varphi_j\}_{j\in J} : F \to G$ gives rise to a 
	unique isomorphism 
	$$\widetilde{\varphi} : \lim\limits_{j\in J} F_j \longrightarrow \lim\limits_{j\in J} G_j$$ 
	commuting the following diagrams for all $j \in J$. 
	\begin{align*}
		\xymatrix{
			F_j \ar[rr] \ar[d]_{\varphi_j} && \lim\limits_{j\in J} F_j \ar[d]^{\widetilde{\varphi}} \\
			G_j \ar[rr] && \lim\limits_{j\in J} G_j
		}
	\end{align*}
	This being a straightforward application of the universal property of colimits, 
	we omit its proof. 
\end{remark}

%%%%%%%%%%%%%%%%%%%%%%%%%%%%%%
\section{Biswas-Borne correspondence for equivariant connections}
In this section, we assume that $D$ is a smooth irreducible Cartier divisor on $X$. 
We now describe the notion of logarithmic connections on $(X, D)$. 
We denote by $\Omega^1_{X/\bb C}(\log D)$ the sheaf of meromorphic differentials on $X$ 
having at most logarithmic poles along $D$.

\subsection{Parabolic connections}
Let $E$ be a vector bundle on $X$. 

\begin{definition}\label{def:logarithmic-connection}
	\begin{enumerate}[(i)]
		\item A \textit{logarithmic connection} on $E$ is a $\bb C$-linear sheaf homomorphism
		\begin{align}
			\nabla : E \to E \otimes_{\mc O_X} \Omega^1_{X/\bb C}(\log D)
		\end{align}
		satisfying the Leibniz rule: 
		$$\nabla(f\cdot s) = f\nabla s + s\otimes df,$$
		for all locally defined sections $f$ and $s$ of $\mc O_X$ and $E$, respectively. 
		It is usually denoted as a pair $(E, \nabla)$. 
		If the image of $\nabla$ lands in the subsheaf $E\otimes\Omega^1_{X/\bb C}$, 
		then $\nabla$ is said to be a {\it holomorphic connection} on $E$. 
		
		\item A \textit{morphism of logarithmic connections} $(E, \nabla)\to (E', \nabla')$ 
		is given by a vector bundle morphism $\phi : E \longrightarrow E'$ such that the 
		following diagram commutes: 
		\begin{equation*}\label{diag:morphism-of-log-connections-defn}
			\begin{gathered}
				\xymatrix{
					E \ar[rr]^-{\nabla} \ar[d]_{\phi} && E\otimes\Omega^1_{X/\bb C}(\log D) 
					\ar[d]^{\phi\otimes\Id} \\
					E' \ar[rr]^-{\nabla'} && E'\otimes\Omega^1_{X/\bb C}(\log D)\,. 
				}
			\end{gathered}
		\end{equation*}
	\end{enumerate}
\end{definition}

Let ${\Con}(X, D)$ be the category of logarithmic connections on $(X, D)$.

\begin{definition}\label{def:parabolic-connection}
	A \textit{parabolic connection} on $(X, D)$ having rational parabolic weights 
	in $\frac{1}{r}\bb Z$ is a functor 
	$$(E_\bullet, \nabla_\bullet) : \left(\frac{1}{r}\bb Z\right)^{\rm op} 
	\longrightarrow {\Con}(X, D)$$ 
	such that for each $n \in \bb N$ there is a natural isomorphism of functors 
	$$\left(E_\bullet\otimes_{\mc O_X}{\mc O}_X(-nD)\ ,\  \nabla_\bullet(-nD)\right)\overset{\sim}{\rightarrow} \left(E_{\bullet+n}\ ,\  \nabla_{\bullet+n}\right)$$
	such that the following diagram commutes 
	\begin{equation}
		\begin{gathered}
			\xymatrix{
				\left(E_{\bullet+n}\,,\,\nabla_{\bullet+n}\right) \ar[rr] && \left(E_{\bullet}\ ,\ \nabla_{\bullet}\right) \\
				\left(E_{\bullet}\otimes_{\mc{O}_X}\mc{O}_X(-nD)\ ,\ \nabla_{\bullet}(-nD)\right) \ar[u]^{\simeq} \ar[rru]_(0.6){\textnormal{Id}_{E_{\bullet}}\otimes\iota_n} &&
			}
		\end{gathered}
	\end{equation}
	where $\iota_n$ stands for the canonical inclusion $\mc{O}(-nD)\subset \mc{O}_X$.
\end{definition} 

\subsection{Equivariant parabolic connection}%\hfill\\
Let $\sigma:G\times X\rightarrow X$ be an action of a connected affine algebraic group $G$ 
on $X$ and satisfying $\sigma^*(D) = p_2^*(D)$. For each $g \in G$, let 
$$\sigma_g : X \longrightarrow X$$ be the map given by multiplication by $g$. 
Consider the sheaf $\Omega^1_{X/\bb{C}}(\log D)$ on $X$, which is locally free 
\cite[Properties 2.2]{EsnVieh92}. Clearly, the assumption $\sigma^*(D) = p_2^*(D)$ 
gives rise to a linearization on $\Omega^1_{X/\bb{C}}(\log D)$: 
$$\psi : p_2^*\left(\Omega^1_{X/\bb{C}}(\log D)\right) \stackrel{\simeq}{\longrightarrow} 
\sigma^*\left(\Omega^1_{X/\bb{C}}(\log D)\right).$$ 
Thus, for each $g\in G$, we get an isomorphism 
$$\psi_g : \Omega^1_{X/\bb{C}}(\log D) \stackrel{\simeq}{\longrightarrow} 
\sigma_g^*\left(\Omega^1_{X/\bb{C}}(\log D)\right)$$
of locally free sheaves. Moreover, the canonical sheaf morphism 
$$\sigma_g^*\left(\Omega^1_{X/\bb{C}}(\log D)\right) 
\longrightarrow \Omega^1_{X/\bb{C}}(\log \sigma_g^*(D))$$
is an isomorphism for each $g\in G$ \cite[Lemma 3.5]{BorLaa23}. 
Composing these isomorphisms, we get the following collection of 
isomorphisms of locally free sheaves:
\begin{equation}\label{eqn:isomorphism-on-logarithmic-differentials}
	\delta_g : \Omega^1_{X/\bb{C}}(\log D) \stackrel{\simeq}{\longrightarrow} 
	\Omega^1_{X/\bb{C}}(\log \sigma_g^*(D)), \ \forall\ g \in G. 
\end{equation}
Following \cite{BisHe13} we define the notion of $G$-equivaraint logarithmic connection 
and $G$-equivariant parabolic connection as follows.

\begin{definition}\label{def:equivariant-logarithmic-and-parabolic-connection}
	\begin{enumerate}[(i)]%\mbox{}
		\item A \textit{$G$-equivariant logarithmic connection on} $(X, D)$ is given by 
		a logarithmic connection $(E, \widehat{\nabla})$ on $(X, D)$ together with a 
		collection of isomorphisms 
		\begin{align*}
			\phi_g : E \stackrel{\simeq}{\longrightarrow} \sigma_g^*(E), \ \forall\ g\in G, 
		\end{align*}
		making the following diagram commutes for all $g \in G$: 
		\begin{equation}\label{diag:equivariant-logarithmic-connection-diagram}
			\begin{gathered}
				\xymatrix{
					E \ar[rr]^-{\widehat{\nabla}} \ar[d]_{\phi_{g}} && 
					E\otimes_{\mc{O}_X}\Omega^1_{X/\bb{C}}(\log D) 
					\ar[d]^{\phi_{g}\otimes\delta_g} \\
					\sigma_g^*(E) \ar[rr]^-{\sigma_g^*(\widehat{\nabla})} && 
					\sigma_g^*(E)\otimes_{\mc{O}_X}\Omega^1_{X/\bb{C}}(\log \sigma_g^*(D)), 
				}
			\end{gathered}
		\end{equation}
		where $\delta_g$ is the isomorphism in \eqref{eqn:isomorphism-on-logarithmic-differentials}.
		\item A \textit{$G$-equivariant parabolic connection on} $(X, D)$ is given by 
		a parabolic connection $(E_{\bullet}, \widehat{\nabla}_{\bullet})$ on $(X, D)$ 
		(see Definition \ref{def:parabolic-connection}) equipped with a collection of 
		natural isomorphisms 
		\begin{align}
			\phi_{h,\bullet} : E_{\bullet} \simeq \sigma_h^*(E_{\bullet})\ \ \forall \ h\in G
		\end{align}
		satisfying $\sigma_g^*(\phi_{h, \bullet}) \circ \phi_{g, \bullet} = \phi_{hg, \bullet}$ 
		for all $g, h \in G$, where we canonically identify 
		$\sigma_g^*(\sigma_h^*(E_{\bullet})) \simeq \sigma_{hg}^*(E_{\bullet})$. 
		These natural isomorphisms are also required to make the following diagrams commute 
		for all $h \in G$: 
		\begin{align}\label{diag:equivariant-parabolic-connection-diagram}
			\begin{gathered}
				\xymatrix{
					E_{\bullet} \ar[rr]^-{\widehat{\nabla}_{\bullet}} \ar[d]_{\phi_{h, \bullet}} 
					&& E_{\bullet}\otimes_{\mc{O}_X}\Omega^1_{X/\bb{C}}(\log D) 
					\ar[d]^{\phi_{h, \bullet}\otimes\delta_h} \\
					\sigma_h^*(E_{\bullet}) \ar[rr]^-{\sigma_h^*(\widehat{\nabla}_{\bullet})} 
					&& \sigma_h^*(E_{\bullet})\otimes_{\mc{O}_X}\Omega^1_{X/\bb{C}}(\log \sigma_h^*(D)), 
				}
			\end{gathered}
		\end{align}
		where $\delta_h$ is the isomorphism as 
		in \eqref{eqn:isomorphism-on-logarithmic-differentials}.
	\end{enumerate}
\end{definition}

The following result shows existence of a canonical logarithmic connection on the 
line bundle $\mc{O}_X(D)$ which we use later (c.f. \cite{EsnVieh92}). 

\begin{lemma}\cite[Lemma 3.7]{BorLaa23}\label{lem:canonical-logarithmic-connection}
	Let $B = \sum_{i\in I}\mu_i D_i$ be a Cartier divisor with support in $D$. 
	Then there exists a canonical logarithmic connection: 
	\begin{align}
		d(B) : \mc{O}_X(B) \rightarrow \mc{O}_X(B) \otimes_{\mc{O}_X} \Omega^1_{X/\bb C}(\log D)
	\end{align}
	characterized by the equation 
	$d(B)\left(\prod_{i\in I}x_i^{-\mu_i}\right) 
	= -\sum_{i\in I}x_i^{-\mu_i}\cdot\sum_{i\in I}\mu_i\dfrac{dx_i}{x_i}$, 
	where $x_i$ is a local equation of $D_i$. 
\end{lemma}

\begin{lemma}
	With the assumptions as in the beginning of this section, the pair $(\mc O_X(D), d(D))$ 
	is a  $G$-equivariant logarithmic connection on $(X, D)$ in the sense of 
	Definition \ref{def:equivariant-logarithmic-and-parabolic-connection}.
\end{lemma}

\begin{proof}
	The conditions $\sigma_g^*(D) = D$ give rise to a collection of isomorphims 
	$\theta_g : \mc{O}(D) \stackrel{\simeq}{\longrightarrow} \sigma_g^*\mc{O}_X(D)$ 
	of line bundles. These conditions imply that the local equation of $D$ differ 
	by multiplying by a unit after composing with $\sigma_g$. From this, it is easy 
	to see that the diagram \eqref{diag:equivariant-logarithmic-connection-diagram} 
	for $(\mc{O}_X(D), d(D))$ commutes by taking $\phi_g = \theta_g$ and 
	$\widehat{\nabla} = d(D)$ in 
	Definition \ref{def:equivariant-logarithmic-and-parabolic-connection}.
\end{proof}

\begin{remark}\label{rem:twist-of-a-logarithmic-connection}
	There is a well-defined notion of tensor product of two logarithmic connections, 
	which is again a logarithmic connection \cite[3.3.2]{BorLaa23}. 
	Thus, given a vector bundle $E$ admitting a logarithmic connection $\nabla$ with 
	poles along $D$, and a Cartier divisor $B$ with support in $D$, one can {\it twist} 
	the connection $\nabla$ by $B$ by considering the tensor product connection on 
	$E\otimes\mc{O}_X(B)$. We shall denote this twisted connection by $\nabla(B)$. 
	The tensor product of two $G$-equivariant logarithmic connections is again $G$-equivariant.
\end{remark}

\subsection{Equivariant logarithmic connection on root stack}
The tautological invertible sheaf  $\ms{M}$ on the root stack $\mf{X}$ gives rise to a 
smooth irreducible effective Cartier divisor $\mf{D}$ on $\mf{X}$ satisfying 
$\pi^*(D) = r\mf{D}$ (see \cite[Remark 2.7]{BorLaa23}). 
Note that $\ms{M} \simeq \mc{O}_{\mf X}(\mf{D})$. 
As we have seen in Proposition \ref{prop:lifting-action-on-root-stack}, under the 
additional assumption that $G$ has no non-trivial character, the $G$-action can be 
lifted to a $G$--action $\sigma_{\mf X} : G\times\mf{X} \longrightarrow \mf{X}$ 
such that the natural map $\pi : \mf{X} \longrightarrow X$ is $G$-equivariant. 
The notion of logarithmic connections with poles along a strict normal crossings divisor makes 
sense for smooth Deligne-Mumford stacks as well. Thus, we can speak about equivariant logarithmic 
connections on the root stack $\mf{X}$. Now, as in the case of schemes, the linearization on 
$\mc{O}_{\mf X}(\mf{D})$ gives rise to the following isomorphisms on the sheaf of 
logarithmic differentials on the root stack: 
\begin{equation}\label{eqn:stacky-logarithmic-differential-isomorphism}
	\eta_g : \Omega^1_{\mf{X}/\bb{C}}(\log \mf{D}) \stackrel{\simeq}{\longrightarrow} 
	\Omega^1_{\mf{X}/\bb{C}}(\log \sigma_{\mf{X}, g}^*(\mf{D})), \ \forall\ g \in G. 
\end{equation}

\begin{definition}\label{def:equivariant-stacky-logarithmic-connection}
	A \textit{$G$-equivariant logarithmic connection on} $(\mf{X}, \mf{D})$ is a 
	logarithmic connection $$\nabla : \mc{E} \longrightarrow 
	\mc{E}\otimes_{{\mc O}_{\mf X}}\Omega^1_{\mf{X}/\bb{C}}(\log\ \mf{D})$$
	such that $\mc{E}$ is a $G$-equivariant vector bundle on $\mf{X}$ 
	admitting vector bundle isomorphisms 
	$$\varphi_g : \mc{E} \stackrel{\simeq}{\longrightarrow} \sigma_{\mf X,g}^*(\mc E), 
	\ \forall\  g\in G,$$
	such that the following diagram commutes for all $g \in G$: 
	\begin{equation}\label{diag:G-connection-diagram}
		\begin{gathered}
			\xymatrix{
				\mc{E}\ar[rr]^-{\nabla}\ar[d]_{\varphi_g} && 
				\mc{E}\otimes_{{\mc O}_{\mf X}}\Omega^1_{\mf{X}/\bb{C}}(\log \mf{D}) 
				\ar[d]^{\varphi_g\otimes\eta_g} \\
				\sigma_{\mf X, g}^*(\mc{E}) \ar[rr]^-{\sigma_{\mf X, g}^*(\nabla)} && 
				\sigma_{\mf X, g}^*(\mc{E})\otimes_{{\mc O}_{\mf X}} 
				\Omega^1_{\mf{X}/\bb{C}}(\log \sigma_{\mf{X}, g}^*(\mf{D})), 
			}
		\end{gathered}
	\end{equation}
	where $\eta_g$ is the isomorphism defined in 
	\eqref{eqn:stacky-logarithmic-differential-isomorphism}.
	
	We shall denote by $\textbf{Con}^{G}(\mf{X}\ ,\ \mf{D})$ the category of $G$-equivariant logarithmic connections on $(\mf{X}\ ,\ \mf{D})$\ .
	
\end{definition}
\begin{remark}\label{rem:canonical-equivariant-logaritmic-connection-on-stack}
	Proposition \ref{prop:linearization-on-tautol-bun-ovr-rt-stack} give rise to isomorphisms of line bundles $\Theta_g : \mc{O}_{\mf X}(\mf {D})\simeq \sigma_{\mf{X},g}^*\mc{O}_{\mf X}(\mf{D})$ for each $g\in G$. Now, just as in Lemma \ref{lem:canonical-logarithmic-connection}, $\mc{O}_{\mf X}(\mf{D})$ admits a canonical logarithmic connection $d(\mf{D})$, which can be constructed by considering \'etale atlases on $\mf X$. It is not difficult to see that this connection is also $G$-equivariant in the sense of Definition \ref{def:equivariant-stacky-logarithmic-connection} by taking $\varphi_g=\Theta_g$.
\end{remark}

\subsection{Equivariant correspondence}
\label{sec:biswas-borne-correspondence-for-equivariant-connections}

It is shown in \cite{BorLaa23} that there is an equivalence between the category of 
logarithmic connections on $(\mf X, \mf D)$ and the category of parabolic connections 
on $(X, D)$ having weights in $\frac{1}{r}\bb{Z}$, and under this equivalence, 
the category of holomorphic connections on $\mf X$ is equivalent to the category 
of strongly parabolic connections on $(X, D)$. It is given as follows.

\begin{enumerate}[(i)]
	\item \cite[Definition 4.18]{BorLaa23} 
	To each logarithmic connection $(\mc E, \nabla)$ on $(\mf X, \mf D)$, one associates 
	a parabolic connection $(E_{\bullet}, \widehat{\nabla}_{\bullet})$ on $(X, D)$, where 
	\begin{enumerate}[(a)]
		\item $E_{\bullet}$ is the parabolic bundle given as in 
		Theorem \ref{thm:Biswas-Borne-correspondence}, and 
		
		\item $\widehat{\nabla}_{\bullet}$ is given by 
		$$\widehat{\nabla}_{\frac{\ell}{r}} = \pi_*(\nabla(-\ell\mf D)),$$
		where $\nabla(-\ell\mf D)$ is the tensor product of the connections 
		$\nabla$ an $d(-\ell\mf D)$. 
		
		(Note that, the push-forward of connection makes sense because the 
		morphism of log-pairs $(\pi, r) : (\mf X, \mf D) \longrightarrow (X, D)$ 
		is log-\'etale, see \cite[Remark 3.8 and Remark 4.19]{BorLaa23}). 
	\end{enumerate}
	
	\item \cite[Definition 4.29]{BorLaa23} 
	On the other hand, for any such parabolic connection 
	$(E_{\bullet}, \widehat{\nabla}_{\bullet})$ on $(X, D)$, one associates the 
	logarithmic connection $(\mc E, \nabla)$ on $(\mf X, \mf D)$, where 
	$$\mc E = \int^{\frac{1}{r}\bb{Z}} \pi^*(E_{\bullet})\otimes_{\mc{O}_{\mf{X}}} {{\mc O}(\bullet r\mf D)}$$ 
	stands for the coend as in Theorem \ref{thm:Biswas-Borne-correspondence}, 
	and $\nabla$ is the unique connection on the coend $\mc E$ compatible with the 
	given connections at each term of the coend (the existence of such a connection 
	is guaranteed by \cite[Lemma 4.27]{BorLaa23}). 
\end{enumerate}

Now we have the following lemma that will be use in the proof of the next theorem. 

\begin{lemma}\label{lem:logarithmic-differential-isomorphism}
	The induced morphism 
	$\pi^*\Omega^1_{X/\bb C}(log(D)) \rightarrow \Omega^1_{{\mf X}/{\bb C}}(log(\mf D))$ 
	under the log-\'etale morphism $\pi$ is an isomorphism. 
\end{lemma}

\begin{proof}
	This is basically the stacky version of \cite[Lemma 3.5]{BorLaa23}, 
	and a similar proof works by considering \'etale atlases for $\mf X$. 
\end{proof}

\begin{theorem}\label{thm:equivariant-connection-correspondence}
	Let $X$ be an irreducible smooth complex projective variety with an action of a 
	connected affine algebraic group $G$. Let $D$ be a smooth irreducible Cartier 
	divisor invariant under the $G$-action. We give $\mf{X}$ the induced $G$-action 
	from Proposition \ref{prop:lifting-action-on-root-stack}. Then there is an 
	equivalence between the category of $G$-equivariant logarithmic connections 
	on $(\mf{X}, \mf{D})$ and the category of $G$-equivariant parabolic connections 
	on $(X,D)$. 
\end{theorem}

\begin{proof}
	Let $\sigma: G\times X\rightarrow X$ and $\sigma_{\mf X} : G\times \mf{X}\rightarrow\mf{X}$ 
	be the action maps. Let $(\mc{E}, \nabla)$ be a $G$-equivariant logarithmic connection 
	on $(\mf{X},\mf{D})$. This means that we have the following commutative diagram for 
	each $g\in G$ (see Definition \ref{def:equivariant-stacky-logarithmic-connection}):
	\begin{align}
		\xymatrix{
			\mc{E} \ar[rr]^-{\nabla} \ar[d]_{\varphi_g}^{\simeq} && 
			\mc{E}\otimes_{{\mc O}_{\mf X}}\Omega^1_{\mf{X}/\bb{C}}(\log \mf{D}) 
			\ar[d]^{\varphi_g\otimes\eta_g} \\ 
			\sigma_{\mf{X}, g}^*(\mc{E}) \ar[rr]^-{\sigma_{\mf{X}, g}^*(\nabla)} && 
			\sigma_{\mf{X}, g}^*(\mc{E})\otimes_{{\mc O}_{\mf X}} 
			\Omega^1_{\mf{X}/\bb{C}}(\log \sigma_{\mf X,g}^*(\mf{D})) 
		}
	\end{align}
	For each $g \in G$, let 
	$$\Theta_g : \mc{O}_{\mf X}(\mf{D}) \stackrel{\simeq}{\longrightarrow} 
	\sigma_{\mf X,g}^*(\mc{O}_{\mf X}(\mf{D}))$$
	be the isomorphisms resulting from the linearization on $\mc{O}_{\mf X}(\mf{D})$ 
	in Proposition \ref{prop:linearization-on-tautol-bun-ovr-rt-stack}. 
	The canonical logarithmic connection $(\mc{O}_{\mf X}(\mf{D}), d(\mf D))$ is 
	$G$-equivariant. For each integer $\ell$, we can take the tensor product of $\nabla$ 
	with the canonical connection $d(-\ell\mf{D})$ on $\mc{O}_{\mf X}(-\ell\mf{D})$, 
	which leads us to the following commutative diagram for all $g \in G$:
	\begin{align}\label{diag:diagram-3}
		\xymatrix{
			\mc{E}(-\ell\mf{D}) \ar[rrr]^-{\nabla(-\ell\mf{D})} 
			\ar[d]_{\varphi_g\otimes \Theta_g^{-\ell}}^{\simeq} &&& 
			\mc{E}(-\ell\mf{D}) \underset{\mc{O}_{\mf X}}{\bigotimes} 
			\Omega^1_{\mf{X}/\bb{C}}(\log \mf{D}) 
			\ar[d]^{(\varphi_g\otimes\Theta_g^{-\ell})\otimes\eta_g} \\ 
			\sigma_{\mf X, g}^*(\mc{E}(-\ell\mf{D})) 
			\ar[rrr]^-{\sigma_{\mf X, g}^*(\nabla(-\ell\mf{D}))} &&& 
			\sigma_{\mf X, g}^*(\mc{E}(-\ell\mf{D})) \underset{\mc{O}_{\mf X}}{\bigotimes} 
			\Omega^1_{\mf{X}/\bb{C}}(\log \sigma_{\mf{X}, g}^*(\mf{D})) 
		}
	\end{align}
	Now, recall the proof of Theorem \ref{thm:equiv-Biswas-Borne-correspondence}, 
	where we have seen the isomorphism of functors 
	$$\sigma^*\circ\pi_*\simeq(\Id\times\pi)_*\circ\sigma_{\mf X}^*.$$ 
	This gives rise to isomorphisms $\sigma_g^*\circ\pi_*\simeq \pi_*\circ\sigma_{\mf X,g}^*$ 
	for all $g \in G$. Also, by Lemma \ref{lem:logarithmic-differential-isomorphism} 
	we have the isomorphisms: 
	\begin{align*}
		& \pi^*\Omega^1_{X/\bb{C}}(\log D) \stackrel{\simeq}{\longrightarrow} 
		\Omega^1_{\mf{X}/\bb{C}}(\log \mf{D}), \\ 
		\text{and } \ \ 
		& \pi^*\Omega^1_{X/\bb{C}}(\log \sigma_g^*(D)) \stackrel{\simeq}{\longrightarrow}
		\Omega^1_{\mf{X}/\bb{C}}(\log\ \sigma_{\mf{X},g}^*(\mf{D})). 
	\end{align*}
	Applying $\pi_*$ on the diagram \eqref{diag:diagram-3} and using the projection formula, 
	we have the following commutative diagram for each $g \in G$: 
	\begin{align}
		\xymatrix{
			E_{\frac{\ell}{r}} \ar[rr]^-{\widehat{\nabla}_{\frac{\ell}{r}}} 
			\ar[d]_{\phi_{g, \frac{\ell}{r}}}^{\simeq} && E_{\frac{\ell}{r}} 
			\otimes_{{\mc O}_{X}} \Omega^1_{X/\bb{C}}(\log D) 
			\ar[d]^{(\phi_{g, \frac{\ell}{r}})\otimes \Id} \\ 
			\sigma_g^*(E_{\frac{\ell}{r}})\ar[rr]^-{\sigma_g^*(\widehat{\nabla}_{\frac{\ell}{r}})} 
			&& \sigma_g^*(E_{\frac{\ell}{r}}) \otimes_{{\mc O}_{X}} \Omega^1_{X/\bb{C}}(\log D)
		}
	\end{align} 
	This produces a $G$-equivariant connection $(\widehat{\nabla}_{\bullet})$ on the 
	resulting $G$-equivariant parabolic bundle $E_{\bullet}$.
	
	Conversely, assume that $(E_{\bullet}, \widehat{\nabla}_{\bullet})$ is a 
	$G$-equivariant parabolic connection on $(X, D)$. We have already seen in 
	Theorem \ref{thm:equiv-Biswas-Borne-correspondence} that the coend 
	$\mc{E} = \int^{\frac{1}{r}\bb{Z}}\pi^*(E_{\bullet})\otimes\mc{O}_{\mf X}(\bullet\,r\mf{D})$ 
	admits a $G$-equivariance structure, and moreover, by \cite[Lemma 4.27]{BorLaa23} it admits 
	a unique logarithmic connection $\nabla$ compatible with the given logarithmic connections 
	$\pi^*(\widehat{\nabla}_{\bullet})\otimes d(\bullet r\mf{D})$ on the bundles 
	$\pi^*(E_{\bullet})\otimes_{\mc{O}_{\mf{X}}}\mc{O}_{\mf X}(\bullet r\mf{D})$. 
	In fact, it follows from the proof of \cite[Lemma 4.27]{BorLaa23} that 
	$(\mc{E}, \nabla)$ is the coend of the bifunctor 
	$$\left(\pi^*(E_{\bullet})\otimes\mc{O}_{\mf X}(\bullet r\mf{D}), 
	\pi^*(\widehat{\nabla}_{\bullet})\otimes d(\bullet\,r\mf{D})\right) : 
	\left(\frac{1}{r}\bb{Z}\right)^{op}\times \left(\frac{1}{r}\bb{Z}\right) 
	\longrightarrow {\Con}(\mf{X}, \mf{D}).$$ 
	Similarly, for each $g \in G$ we can interpret the connection 
	$\left(\sigma_{\mf X, g}^*(\mc{E}), \sigma_{\mf X, g}^*(\nabla)\right)$ 
	as a coend of the bifunctor 
	$$\left(\sigma_{\mf{X}, g}^*(\pi^*(E_{\bullet})\otimes\mc{O}_{\mf X}(\bullet\,r\mf{D})), 
	\sigma_{\mf{X}, g}^*(\pi^*(\widehat{\nabla}_{\bullet})\otimes d(\bullet\,r\mf{D}))\right) 
	: \left(\frac{1}{r}\bb{Z}\right)^{op}\times \left(\frac{1}{r}\bb{Z}\right) 
	\longrightarrow {\Con}({\mf X}, \mf{D}).$$ 
	Since $(E_{\bullet}, \widehat{\nabla}_{\bullet})$ is $G$--equivariant, we have 
	commutative diagrams \eqref{diag:equivariant-parabolic-connection-diagram} 
	for each $h\in G$. Applying $\pi^*$ to these diagrams, and using 
	$\pi^*\circ\sigma_h^* = \sigma_{\mf{X}, h}^*\circ\pi^*$ together with the 
	Lemma \ref{lem:logarithmic-differential-isomorphism}, 
	we get a commutative diagram of the form
	\begin{equation}\label{diag:diagram-4}
		\begin{gathered}
		\xymatrix{
			\pi^*(E_{\bullet}) \ar[rrr]^-{\pi^*(\widehat{\nabla}_{\bullet})} 
			\ar[d]_{\pi^*(\phi_{h,\bullet})} &&& 
			\pi^*(E_{\bullet})\ \underset{\mc{O}_{\mf X}}{\bigotimes} 
			\Omega^1_{\mf{X}/\bb{C}}(\log\ \mf{D}) \ar[d]^{\pi^*(\phi_{h,\bullet})\otimes\eta_h} \\
			\sigma_{\mf{X}, h}^*(E_{\bullet}) \ar[rrr]^-{\pi^*(\sigma_h^*(\widehat{\nabla}_{\bullet}))} 
			&&& \sigma_{\mf{X}, h}^*(\pi^*(E_{\bullet})) \underset{\mc{O}_{\mf X}}{\bigotimes} 
			\Omega^1_{\mf{X}/\bb{C}}(\log\ \sigma_{\mf{X}, h}^*(\mf{D}))
			}
		\end{gathered}
	\end{equation}
	where $\eta_h$ is as in \eqref{eqn:stacky-logarithmic-differential-isomorphism}.
	We now consider the tensor product of the two horizontal arrows in the diagram \eqref{diag:diagram-4} with the functor 
	$$(\mc{O}_{\mf X}(\bullet\,r\mf{D}), d(\bullet\,r\mf{D})) : 
	\left(\frac{1}{r}\bb{Z}\right)^{\rm op} \longrightarrow {\Con}^G({\mf X}, {\mf D}),$$ 
	which will give rise to diagrams of the form 
	\begin{equation*}
		\xymatrix{
			\pi^*(E_{\bullet})\ \underset{\mc{O}_{\mf X}}{\bigotimes} \mc{O}_{\mf{X}}(\bullet r\mf{D}) 
			\ar[rrr]^-{\pi^*(\widehat{\nabla}_{\bullet})\otimes d(\bullet\ r\mf{D})} 
			\ar[d]_{\pi^*(\phi_{h, \bullet})\otimes\Theta_{h}^{\otimes\bullet}} &&& 
			\pi^*(E_{\bullet}) \underset{\mc{O}_{\mf X}}{\bigotimes} \mc{O}_{\mf X}(\bullet r\mf{D}) 
			\underset{\mc{O}_{\mf X}}{\bigotimes} \Omega^1_{\mf{X}/\bb{C}}(\log \mf{D}) 
			\ar[d]^{(\pi^*(\phi_{h,\bullet}) \otimes \Theta_h^{\otimes\bullet})\otimes\eta_h} \\
			\sigma_{\mf{X},h}^*(E_{\bullet})\ \underset{\mc{O}_{\mf X}}{\bigotimes} 
			\sigma_{\mf{X},h}^*(\mc{O}_{\mf X}(\bullet\ r\mf{D})) 
			\ar[rrr]^-{\pi^*(\sigma_h^*(\widehat{\nabla}_{\bullet}))\otimes d(\bullet\,r\mf{D})} 
			&&& \sigma_{\mf{X}, h}^*(\pi^*(E_{\bullet})) \underset{\mc{O}_{\mf X}}{\bigotimes} 
			\mc{O}_{\mf X}(\bullet\,r\mf{D}) \underset{\mc{O}_{\mf X}}{\bigotimes}
			\Omega^1_{\mf{X}/\bb{C}}(\log \sigma_{\mf{X}, h}^*(\mf{D})) 
		}
	\end{equation*}
	where $\Theta_h$ is as in \eqref{rem:canonical-equivariant-logaritmic-connection-on-stack}. 
	Each term in the diagram above is a bifunctor 
	$$\left(\frac{1}{r}\bb{Z}\right)^{\rm op}\times \left(\frac{1}{r}\bb{Z}\right) 
	\longrightarrow \textbf{Con}(\mf{X}, \mf{D}).$$
	Finally, taking the coends of these bifunctors, we get a commutative diagram of the form 
	\begin{align*}
		\xymatrix{
			\mc{E} \ar[rr]^(.34){\nabla} \ar[d]_{\varphi_h} && 
			\mc{E}\otimes_{{\mc O}_{\mf X}}\Omega^1_{\mf{X}/\bb{C}}(\log\ \mf{D}) 
			\ar[d]^{\varphi_h\otimes\eta_h}\\
			\sigma_{\mf{X},h}^*(\mc{E}) \ar[rr]^(.3){\sigma_{\mf{X},h}^*(\nabla)} 
			&& \sigma_{\mf{X},h}^*(\mc{E})\otimes_{\mc{O}_{\mf{X}}} 
			\Omega^1_{\mf{X}/\bb{C}}(\log \sigma_{\mf{X},h}^*(\mf{D}))
		}
	\end{align*} 
	Thus, under the correspondence of \cite{BorLaa23} (see the beginning of 
	\S\,\ref{sec:biswas-borne-correspondence-for-equivariant-connections}), 
	the resulting logarithmic connection $(\mc{E}, \mf{D})$ is $G$-equivariant. 
	This completes the proof. 
\end{proof}

\begin{remark}
	For $\lambda \in \bb{C}$, there is a notion of {\it logarithmic $\lambda$--connections} 
	on $X$, which are pairs $(E, \nabla)$, where $E$ is a vector bundle on $X$ and 
	$\nabla : E \longrightarrow E\otimes\Omega^1_{X/\bb{C}}(\log D)$ is a $\bb{C}$-linear 
	homomorphism of sheaves satisfying the {\it $\lambda$--twisted Leibniz rule}: 
	$$\nabla(fs) = f\nabla(s) + \lambda\cdot df\otimes s,$$
	for all locally defined sections $f$ and $s$ of $\mc O_X$ and $E$, respectively. 
	Such connections can be twisted by a Cartier divisor $B$ with support in $D$. 
	The notion of logarithmic $\lambda$--connection gives rise to the notion of 
	{\it parabolic $\lambda$--connection} generalizing the 
	Definition \ref{def:equivariant-logarithmic-and-parabolic-connection}. 
	When $X$ is endowed with an action of a connected affine algebraic group $G$ 
	having no non-trivial characters, and $D$ is invariant under the given $G$--action 
	on $X$, we may talk about $G$--equivariant parabolic $\lambda$--connection as well. 
	All these notions have their appropriate counterparts on the associated $r$-th root stack 
	$\mf X = \mf X_{(\mc O_X(D), s_D, r)}$. 
	Then Theorem \ref{thm:equivariant-connection-correspondence} can be generalized in the 
	setup of $G$--equivariant $\lambda$--connections following the same line of arguments. 
	The case of $\lambda = 0$ is of particular interest, as it leads to the notion of 
	{\it Higgs bundles}. One may use this to generalize the Biswas-Borne correspondence 
	for Higgs bundles (see \cite{BisMajWon13}) in $G$--equivariant setup. 
\end{remark}

\section*{Acknowledgment}
The authors would like to thank Prof. Indranil Biswas for his helpful comments. The second named author is partially supported by the DST INSPIRE Faculty Fellowship 
(Research Grant No.: DST/INSPIRE/04/2020/000649), Ministry of Science \& Technology, 
Government of India.

%\bibliographystyle{halpha}
%\bibliography{references}

\end{document}